\documentclass[12pt]{article}

\usepackage{a4wide}
\usepackage{amssymb}
\usepackage{amsfonts}
\usepackage{amsmath}
\input xy
\xyoption{arrow} \xyoption{matrix}

\date{}

\newtheorem{proposition}{Proposition}[section]
\newtheorem{theorem}[proposition]{Theorem}
\newtheorem{lemma}[proposition]{Lemma}

\newtheorem{corollary}[proposition]{Corollary}

\def\der{\partial }

\def\nFM0{{\nu }_{F,M_0}}
\def\nFN0{{\nu }_{F,N_0}}
\def\nGN0{{\nu }_{G,N_0}}

\def\N0{ {\bf N}_0 }

\def\t{\otimes}
\def\g{\gamma}

\def\ra{\rightarrow}

\def\Xpm{X^{\pm }}

\def\s{\sigma}
\def\Z{\mathbb{Z}}

\def\l1{{\lambda}_1}

\def\a{\alpha}
\def\a0{ {\alpha }_0}
\def\a1{ {\alpha }_1}

\def\l{\lambda}

\def\nFGM0{{\nu }_{F,G,M_0}}

\def\nFN0{{\nu}_{F,N_0}}

\def\sm{{\sigma}^m}

\def\sm1{{\sigma}^{-1}}

\def\smtp1{{\sigma}^{-t+1}}

\def\S1{S^{-1}}

\def\Xpm1{X^{\pm 1}_1}

\def\sPM1{{\sigma }^{\pm 1}}
\def\sMP1{{\sigma }^{\mp 1 }}

\def\d{\delta}

\def\di{{\rm d.ind}}

\def\L{\Lambda}

\def\CD{{\cal D}}

\def\Ytm1{Y^{t-1}}
\def\Yim1{Y^{i-1}}

\def\CM{{\cal M}}

\def\CF{{\cal F}}

\def\Aut{{\rm Aut}}

\def\Der{{\rm Der }}
\def\ad{{\rm ad }}

\def\ker{ {\rm ker } }

\def\D{ \Delta }

\def\xx{ {\bf x} }
\def\yy{ {\bf y} }

\def\SL2Z{ {\rm SL}_2({\bf Z}) }

\def\Gp1{ G^{1 , 1 } }
\def\P11{ P^{-1 , 1 } }
\def\Pp1{ P^{1 , 1 } }

\def\nCLsr{{}^\nu\kern-2pt {\cal L}^{\sigma , \rho  }}
\def\nP{{}^\nu \kern-2pt P}
\def\nL{{}^\nu\kern-2pt L}
\def\nLL{{}^\nu\kern-2pt \Lambda}
\def\nPsr{{}^\nu\kern-2pt P^{\sigma , \rho  }}
\def\nLsr{{}^\nu\kern-2pt L^{\sigma , \rho  }}
\def\nuCL{{}^\nu\kern-2pt  {\cal L}}
\def\nCLsr{{}^\nu\kern-2pt {\cal L}^{\sigma , \rho  }}
\def\nCL1m{{}^\nu\kern-2pt {\cal L}^{-1 , 1  }}
\def\x1nu{x^\frac{1}{\nu}}
\def\xm1nu{x^{-\frac{1}{\nu}}}

\def\ra{\rightarrow }

\def\nAM0{{\nu }_{{\cal A},M_0}}
\def\nAN0{{\nu }_{{\cal A},N_0}}

\def\Der{ {\rm Der }}

\def\ad{ {\rm ad }}

\def\gm{\mathfrak{m}}

\def\derij{\partial_{{\bf i}, {\bf j}}}

\def\SL{{\rm SL}}

\def\di!{\frac{\der^i}{i!}}
\def\dik!{\frac{\der^k_i}{k!}}

\def\Fp{\mathbb{F}_p}

\def\derik{\der_i^{[k]}}
\def\derij{\der_i^{[j]}}

\def\deril{\der_i^{[l]}}
\def\derjl{\der_j^{[l]}}
\def\derba{\der^{[\alpha ]}}
\def\derbb{\der^{[\beta ]}}

\def\Nn{\mathbb{N}^n}

\def\id{{\rm id}}

\def\ID{{\rm ID}}

\def\N{\mathbb{N}}

\def\0{\overline{0}}
\def\1{\overline{1}}

\def\Ln1{\L_{n,\overline{1}}}

\def\a1{a_{\overline{1}}}

\def\S{\Sigma}

\def\vn1{\overrightarrow{n-1}}

\def\hx{\widehat{x}}

\def\deripk{\der_i^{[p^k]}}

\def\derjpl{\der_j^{[p^l]}}
\def\derijkpk{\der_i^{[j_kp^k]}}
\def\deripkjk{\der_i^{[p^k]j_k}}
\def\CDPn{\CD (P_n)}

\def\FrobCDPn{{\rm Frob}(\CD (P_n))}

\def\FrobsCDPnDn{{\rm Frob}_{s}(\CD (P_n), P_n, \D_n)}
\def\FrobCDPnPnDn{{\rm Frob}(\CD (P_n), P_n, \D_n)}

\def\FrobsCDPnPn{{\rm Frob}_{s}(\CD (P_n), P_n))}
\def\Frob{{\rm Frob}}

\def\ID{{\rm ID}}

\def\derips{\der_i^{[p^s]}}

\def\FPolCDPn{{\rm FPol}(\CDPn )}

\begin{document}

\author{V. V. \  Bavula 
}

\title{Extensions of the Frobenius to the ring of differential operators
on a polynomial algebra in prime characteristic }

\maketitle
\begin{abstract}
Let $K$ be a field of characteristic $p>0$.  It is proved that
each automorphism  $\s \in \Aut_K(\CDPn)$ of the ring $\CDPn$ of
differential operators on a polynomial algebra $P_n= K[x_1, \ldots
, x_n]$ is {\em uniquely} determined by the elements $\s (x_1),
\ldots ,\s (x_n)$;  and  that the set $\Frob (\CDPn)$ of all the
extensions of the Frobenius (homomorphism)  from certain maximal
commutative polynomial subalgebras of $\CDPn$, like $P_n$, to the
ring $\CDPn$ is equal to $\Aut_K(\CDPn ) \cdot \CF$ where $\CF$ is
the set of all the extensions of the Frobenius from $P_n$ to
$\CDPn$ that leave invariant the subalgebra of scalar differential
operators. The set
 $\CF$ is found explicitly, it is large (a typical extension
depends on {\em countably} many independent parameters).

$\noindent $



$${\bf Contents}$$
\begin{enumerate}
\item Introduction. \item Existence and uniqueness of  iterative
$\d$-descent. \item Rigidity of the group $\Aut_K(\CDPn )$. \item
Extensions of the Frobenius to the ring $\CDPn$. \item The sets
$\FrobsCDPnDn $, $s\geq 1$. \item Appendix: Some technical
results.
\end{enumerate}
\end{abstract}


\section{Introduction}
Throughout, ring means an associative ring with $1$, $\N :=\{0, 1,
\ldots \}$ is the set of natural numbers, $p$ is a prime number,
$\Fp := \Z / \Z p$ is the field that contains $p$ elements, $K$ is
an arbitrary field of characteristic $p>0$ (if it is not stated
otherwise), $P_n:= K[x_1, \ldots , x_n]$ is a polynomial algebra,
$\CDPn =\oplus_{\alpha \in \N^n}P_n\der^{[\alpha]}$ is the ring of
differential operators on $P_n$ where
$\der^{[\alpha]}:=\prod_{i=1}^n
\frac{\der_i^{\alpha_i}}{\alpha_i!}$, and  $\D_n=\oplus_{\alpha
\in \N^n}K\der^{[\alpha]}$ is the algebra of {\em scalar}
differential operators on $P_n$.

$\noindent $

{\bf Rigidity of the group of automorphisms $\Aut_K(\CDPn )$}. In
characteristic zero, there is a strong connection between the
groups $\Aut_K(P_n)$ and $\Aut_K(\CDPn )$ as, for example, the
(essential) equivalence of the Jacobian Conjecture for $P_n$ and
the Dixmier Problem/Conjecture for $\CDPn$ shows (see \cite{BCW},
\cite{Tsuchi05}, \cite{Bel-Kon05JCDP}, see also  \cite{JC-DP}).
Moreover, in the class of all the associative algebras conjecture
like the two mentioned conjectures makes sense only for the
algebras $P_m\t \CDPn$ as was proved in \cite{Bav-cpinv} (the two
conjectures can be reformulated in terms of locally nilpotent
derivations that satisfy certain conditions, and the algebras
$P_m\t \CDPn$  are the only associative algebras that have such
derivations). This general conjecture is true iff either the JC or
the DC is true, see \cite{Bav-cpinv}.

In prime characteristic,  relations between the two groups,
$\Aut_K(P_n)$ and $\Aut_K(\CDPn )$, are  even tighter as the
following  result shows.

\begin{theorem}\label{28Jun7}
{\rm (Rigidity of the group $\Aut_K(\CDPn )$)} Let $K$ be a field
of characteristic $p>0$, and $\s , \tau \in \Aut_K(\CDPn )$. Then
$\s = \tau $ iff $\s (x_1) = \tau (x_1) , \ldots , \s (x_n) = \tau
(x_n)$.
\end{theorem}

{\it Remark}.  Theorem \ref{28Jun7} does not hold in
characteristic zero ({\em eg}, the automorphisms $\s : x\mapsto
x+\der $ , $\der \mapsto \der$, and $\tau = \id_{\CD (P_1)}$ of
the first {\em Weyl} algebra $\CD (P_1)=K\langle x, \der \, | \,
\der x- x\der =1\rangle$ are distinct but $\s (x) = \tau (x)$). In
general, in prime characteristic Theorem \ref{28Jun7} does not
hold for localizations of the polynomial algebra $P_n$, {\em eg},
$K[x_1^{\pm 1}, \ldots , x_n^{\pm 1})$, \cite{Bav-ProcAMS2009}.
Note that the $K$-algebra $\CDPn$ is not finitely generated and
neither left nor right Noetherian.

$\noindent $

{\bf Extensions of the Frobenius to the ring of differential
operators $\CDPn$}. In the paper, we are interested in the
question of extending  the {\em Frobenius homomorphism} (the
Frobenius, for short),
$$F:P_n\ra P_n, \;\;  a\mapsto a^p,$$
 to the ring of differential operators
$\CDPn$. There is the canonical one,
$$ F_x:\CDPn \ra \CDPn , \;\; \der^{[\alpha]}\mapsto
\der^{[p\alpha]}, \;\; \alpha \in \N^n , $$ which is called the
{\em canonical Frobenius} on $\CDPn$ that corresponds to the set
of generators $(x)=(x_1, \ldots , x_n)$ for the polynomial algebra
$P_n$.  As a curious fact, note that, by a trivial reason, the
Frobenius $F$ {\em cannot} be extended from the polynomial algebra
$P_n$ to the Weyl algebra $A_n$. To the contrary, as it is proved
in the paper, the set $\Frob (\CDPn , P_n)$ of all extensions of
the Frobenius $F$ on $P_n$ to the algebra $\CDPn$ is massive (a
typical extension depends on countably many parameters).

A ring homomorphism $F': \CDPn \ra \CDPn$ is called a {\em
Frobenius} if it is an extension of the Frobenius from a certain
polynomial subalgebra $P_n= K[x_1', \ldots , x_n']$ of the ring
$\CDPn$ (which is called a {\em Frobenius polynomial subalgebra})
that satisfies natural $\Aut_K(\CDPn )$-invariant conditions like
the canonical Frobenius $F_x$ does (see Section \ref{ETFRB} for
the details). Let $\Frob (\CDPn )$ be the set of all the Frobenius
homomorphisms on $\CDPn$ and $\FPolCDPn $ be the set of all the
Frobenius polynomial subalgebras of $\CDPn$. It is shown that each
Frobenius polynomial subalgebra is a maximal commutative
 subalgebra of $\CDPn$ (Corollary \ref{a21Mar8}.(1)), and so one
 Frobenius polynomial subalgebra cannot contain properly another
 Frobenius polynomial subalgebra.

$\noindent $

{\it Definition}. Let  $\Frob (\CDPn , P_n, \D_n)$ be the set  of
all the extensions $F'$ of the Frobenius $F$ on the polynomial
algebra $P_n$ to $\CDPn$  such that $F'(\D_n ) \subseteq \D_n$.

$\noindent $

The next theorem describes  the sets $\Frob (\CDPn )$, $\Frob
(\CDPn , P_n)$ and $\FPolCDPn $ up to the action of the groups
$\Aut_K(\CDPn )$ and $\Aut_K(P_n)$.

\begin{theorem}\label{18Mar8}
\begin{enumerate}
\item $\Frob (\CDPn )= \Aut_K(\CDPn ) \Frob (\CDPn , P_n, \D_n)$,
i.e. for each Frobenius $F'$ on $\CDPn$ there exists an
automorphism $\s \in \Aut_K(\CDPn )$ and a Frobenius $F''\in\Frob
(\CDPn , P_n, \D_n)$ such that $F'=\s F''\s^{-1}$. \item $\Frob
(\CDPn , P_n )= \Aut_K(P_n) \Frob (\CDPn , P_n, \D_n)$, i.e. for
each Frobenius $F'\in \Frob (\CDPn , P_n )$  there exists an
automorphism $\s \in \Aut_K(P_n)$ and a Frobenius $F''\in\Frob
(\CDPn , P_n, \D_n)$ such that $F'=\s F''\s^{-1}$. \item
$\FPolCDPn  = \Aut_K(\CDPn ) \cdot P_n\simeq \Aut_K(\CDPn ) /
\Aut_K(P_n):= \{ \s \Aut_K(P_n)\, | \, \s \in \Aut_K(\CDPn )\}$,
i.e. for each Frobenius polynomial subalgebra $P_n'$ of $\CDPn$
there exists an automorphism $\s\in \Aut_K(\CDPn )$ such that $\s
(P_n) = P_n'$ and the automorphism $\s$ is unique up to
$\Aut_K(P_n)$.
\end{enumerate}
\end{theorem}

$\noindent $

The set $\Frob (\CDPn , P_n, \D_n)$ is found explicitly (Theorems
\ref{M5Jul7} and \ref{cM5Jul7}). Each Frobenius $F'\in \FrobCDPn$
is not an $\Fp$-algebra isomorphism (though it is a monomorphism
since the ring  $\CDPn$ is simple). Moreover, the ring  $\CDPn$ is
a left and right free  finitely generated $KF'(\CDPn )$-module of
rank $p^{2n}$ (Corollary \ref{7Apr8}). Theorem \ref{18Mar8} shows
that the problem of finding the groups $\Aut_K(P_n)$ and
$\Aut_K(\CDPn )$ is closely related to the problem of finding all
the extensions of the Frobenius from certain polynomial
subalgebras to the ring $\CDPn$.

$\noindent $

{\bf Iterative $\d$-descents}. The question of finding all the
extensions of the Frobenius (as well as many other difficult
questions like the Jacobian Conjecture or the Dixmier Conjecture)
can be reformulated as a question about iterative $\d$-descents.

$\noindent $

{\it Definition}. Let $A$ be a ring and $\d = (\d_1, \ldots
,\d_n)$ be an $n$-tuple of commuting derivations of the ring $A$.
A multi-sequence $\{ y^{[\alpha]}, \alpha \in \N^n\}$ in $A$ is
called an {\bf iterative $\d$-descent} if
$$ y^{[\alpha]}y^{[\beta]}= {\alpha +\beta \choose \alpha}
y^{[\alpha + \beta ]}\;\; {\rm and}\;\; \d^\alpha (y^{[\beta]}) =
y^{[\beta - \alpha ]}\;\; {\rm for \; all}\;\;  \alpha , \beta \in
\N^n.$$

$\noindent $

{\it Example}. $\{ \der^{[\alpha]}, \alpha \in \N^n\}$ is the
iterative $\d$-descent in the ring $\CDPn$ where $\d_1=-\ad (x_1),
\ldots , \d_n= -\ad (x_n)$.

$\noindent $

This concept first appeared in \cite{Bav-simderharp} where using
it the simple derivations of differentially simple Noetherian
commutative rings were classified. One of the key results of the
paper \cite{Bav-simderharp} is the existence and uniqueness of an
iterative $\d$-descent (Theorem 3.8, \cite{Bav-simderharp}). In
the paper, this result is extended to the case of several
commuting derivations and infinite iterative descents (Theorem
\ref{15Nov7}). This result is used in the proofs of  almost all
main results of the paper.

 In Section \ref{EAUAI}, a theory of iterative
$\d$-descents is developed. These results are used freely in the
paper.


\section{Existence and uniqueness of  iterative
$\d$-descent}\label{EAUAI}

The main result of this section is Theorem \ref{15Nov7} on {\em
existence} and {\em uniqueness} of an {\em iterative} $\d$-{\em
descent}. Necessary and sufficient conditions are given (Corollary
\ref{f14Nov7}) for a multi-sequence to be an iterative
$\d$-descent.

$\noindent $

{\bf Iterative $\d$-descents}. Consider the free abelian group of
rank $n$,  $\Z^n=\oplus_{i=1}^n\Z e_i$,  where the set of elements
$e_1=(1,0, \ldots , 0), \ldots , e_n=(0, \ldots , 0,1)$ is the
standard free basis for $\Z^n$. For each element $\alpha =
(\alpha_1, \ldots , \alpha_n)\in \Z^n$, we have $\alpha =
\sum_{i=1}^n \alpha_ie_i$. For elements  $\alpha , \beta \in \Z^n$,
we write $\alpha \geq \beta$ if $\alpha_1 \geq \beta_1, \ldots ,
\alpha_n\geq \beta_n$; we write $\alpha >\beta$ if $\alpha\geq
\beta$ and $\alpha \neq \beta$. For a pair of integers $i$ and $j$
such that $i\leq j$ (resp. $i<j$), let $[i,j]_{dis}:=\Z \cap [i,j]=
\{ i, i+1, \ldots , j\}$ (resp. $[i,j)_{dis}:=\Z \cap [i,j)= \{ i,
i+1, \ldots , j-1\}$).

$\noindent $

{\it Definition}. Let $A$ be a ring and let $\d = (\d_1, \ldots ,
\d_n)$ be an $n$-tuple of commuting derivations of the ring $A$.
For each element $\alpha \in \N^n$, let $\d^\alpha :=
\d_1^{\alpha_1}\cdots \d_n^{\alpha_n}$. A multi-sequence of
elements of the ring $A$, $\yy := \{ y^{[\alpha ] }, \alpha \in
I:=\prod_{i=1}^n[0,l_i)_{dis}\}$, $y^{[0]}:=1$, where $l_i\in \N
\cup \{ \infty \}$, is called a $\d$-{\em descent}  if
$$\d^\alpha (y^{[\beta]})=y^{[\beta -\alpha]} \;\; {\rm for \; all}\;\;  \alpha \in \N^n,
\;\; \beta \in I, $$ where $y^{[\g ]}:=0$ for all $\g \in
\Z^n\backslash \N^n$.

$\noindent $

{\it Definition}. Suppose that $A$ is an $\Fp$-algebra. A
multi-sequence $\{ y^{[\alpha]}, \alpha \in
I:=\prod_{i=1}^n[0,p^{d_i})\}$ where $d_i\in \N \cup\{\infty \}$
is called an {\em iterative} multi-sequence if
$$ y^{[\alpha]}y^{[\beta]}= {\alpha +\beta \choose \beta}y^{[\alpha
+\beta]} \;\; {\rm for \; all}\;\; \alpha , \beta \in I,$$ where
${\alpha +\beta \choose \beta}:= \prod_{i=1}^n {\alpha_i +\beta_i
\choose \beta_i}$ are the multi-binomial coefficients. In prime
characteristic, we view the multi-binomial coefficients as
elements of the field $\Fp$.

$\noindent $

{\it Definition}, \cite{Bav-simderharp}. Let $A$ be an
$\Fp$-algebra. An iterative multi-sequence $\{ y^{[\alpha]},
\alpha \in I=\prod_{i=1}^n[0,p^{d_i})\}$,  $d_i\in \N \cup\{\infty
\}$ in $A$ which is a $\d$-descent is called an {\bf iterative}
$\d$-{\bf descent} of {\em rank} $n$ and  {\em exponent} $(d_1,
\ldots , d_n)$.

$\noindent $

{\it Example}. Let $K$ be a field of characteristic $p>0$, $\CDPn
=\oplus_{\alpha , \beta \in \N^n}Kx^\alpha \derbb$ be the ring of
differential operators on the polynomial algebra $P_n=K[x_1,
\ldots , x_n]$. Consider the inner derivations $\d_1:= - \ad
(x_1), \ldots , \d_n :=-\ad (x_n)$ of the ring $\CDPn$ ($\ad
(x)(y):=xy-yx$) and let $\d : = (\d_1, \ldots , \d_n)$. Then the
multi-sequence $\{ \der^{[\alpha]} :=
\frac{\der_1^{\alpha_1}}{\alpha_1!}\cdots
\frac{\der_n^{\alpha_n}}{\alpha_n!}, \alpha \in \N^n\}$ is an
iterative $\d$-descent (of rank $n$ and  exponent $(\infty ,
\ldots , \infty )$).

$\noindent $

Note that any truncation $\{ y^{[\alpha]}, \alpha \in
I'=\prod_{i=1}^n[0,p^{d_i'})\}$,  $d_i'\leq d_i$, of the iterative
 $\d$-descent $\{ y^{[\alpha]}, \alpha \in I\}$ is  an iterative $\d$-descent of rank $n$ and of exponent $(d_1', \ldots ,
d_n')$.

$\noindent $

{\bf The nil ring of commuting derivations}. A ring $S$ is a {\em
positively filtered} ring if $S$ is a union of its abelian
subgroups, $S= \cup_{i\geq 0} S_i$,  such that $S_0\subseteq
S_1\subseteq \cdots $ and $S_iS_j \subseteq S_{i+j}$ for all
$i,j\geq 0$.  Let $A$ be a ring   and let $\d $ be a
 derivation of the ring $A$. Recall that, for any elements $a,b\in A$
 and a natural number $n$, we have the equality
 $$ \d^n(ab)=\sum_{i=0}^n\, {n\choose i}\d^i(a)\d^{n-i}(b).$$
 It follows directly from this equality that the union of the
abelian groups $N:=N(\d ,A)=\cup_{i\geq 0}\, N_i$, $N_i:= \ker \,
\d^{i+1}$, is a positively
  filtered ring ($N_iN_j\subseteq N_{i+j}$ for
 all $i,j\geq 0$), so-called, the {\em nil ring} of $\d $.
 Clearly,  $N_0= A^\d:=\ker \, \d $  is the
 subring (of {\em constants} for $\d $) of $A$,
   and $N=\{ a\in A \, | \ \d^n (a)=0$
  for some natural $n=n(a)\}$. For all natural  numbers $i$ and $j$,
  $\d^i(N_j) \subseteq N_{j-i}$ where $N_k:=0$ for all $\Z\backslash
  \N$. In general, little is known about  rings $N(\d , A)$, these
  rings have complicated structure and they are not easy objects to
  deal with. Even for the first Weyl algebra there are old open
  problems about these rings, see the paper of J. Dixmier,
  \cite{Dix-1968}. A derivation $\d$ of a ring $A$ is called a {\em
  locally nilpotent} derivation if $A= N(\d , A)$ or, equivalently,
  for each element $a$ of $A$, $\d^i(a)=0$ for all $i\gg 0$. In this
  case, the ring $A=\cup_{i\geq 0}N(\d , A)_i$ is a positively
  filtered ring.

Let $\d = (\d_1, \ldots , \d_n)$ be an $n$-tuple of commuting
derivations of a ring $A$. The intersection $N(\d , A) :=
\cap_{i=1}^n N(\d_i, A)$ is called the {\em nil ring} of $\d$. It
is an $\N^n$-filtered ring, $N(\d , A) = \cup_{\alpha\in \N^n}
N(\d , A, \alpha )$ where
$$ N_\alpha := N(\d , A , \alpha ) := \cap_{i=1}^n N(\d_i, A,
\alpha_i)\;\;\; (N_\alpha N_\beta \subseteq N_{\alpha+\beta} \;\;
{\rm for \; all}\;\;  \alpha, \beta \in \N^n ).$$ Note that $N_0:=
A^\d : = \cap_{i=1}^n A^{\d_i}$ is the subring of $A$, so-called,
the ring of $\d$-constants for $\d$. For $\alpha  ,\beta \in \N^n$,
$\d^\alpha (N_\beta ) \subseteq N_{\beta - \alpha }$ where $N_\g
:=0$ for all $\g \in \Z^n\backslash \N^n$.

Let $\d = (\d_1, \ldots , \d_n)$ be an $n$-tuple of commuting,
locally nilpotent derivations of a ring $A$. Then $A= N(\d , A)$,
and the ring $A= \cup_{\alpha \in \N^n} N(\d, A , \alpha )$ is an
$\N^n$-filtered ring.

$\noindent $

{\it Example}. Let $K$ be a field of characteristic $p>0$, and
$P_n= K[x_1, \ldots , x_n]$ be a polynomial algebra over $K$. Then
$\d := (- \ad (x_1), \ldots , -\ad (x_n))$  is the $n$-tuple of
commuting, locally nilpotent derivations of the ring $\CDPn
=\oplus_{\alpha \in \N^n}P_n\der^{[\alpha]}$, and so the ring
$\CDPn = \cup_{\alpha \in \N^n} N_\alpha$ is $\N^n$-filtered where
$N_\alpha := \oplus_{0\leq \beta \leq \alpha}P_n\der^{[\alpha]}$.

$\noindent $

The following lemma establishes a relation between the filtration
$\{ N_\alpha\}_{\alpha \in \N^n}$ of the nil ring $N(\d , A)$ of $\d
= (\d_1, \ldots , \d_n)$ and $\d$-descents.

\begin{lemma}\label{a13Nov7}
Let $\d = (\d_1, \ldots , \d_n)$ be an $n$-tuple of commuting
derivations of a ring $A$, and $\{ x^{[\alpha]}, \alpha\in I:=
\prod_{i=1}^n [0,d_i]_{dis}\}$ be a $\d$-descent where $d_i\in \N
\cup \{ \infty \}$. Then, for each $\alpha \in I$, $N_\alpha
=\oplus_{0\leq \beta \leq \alpha} A^\d x^{[\beta]}= \oplus_{0\leq
\beta \leq \alpha}x^{[\beta]}A^\d$.
\end{lemma}

{\it Proof}. Let us prove the first equality. For each element
$\alpha \in I$, let $N_\alpha':=\oplus_{0\leq \beta \leq \alpha}
A^\d x^{[\beta]}$. The inclusion $N_\alpha'\subseteq N_\alpha$ is
obvious. To prove the reverse inclusion we use induction on $|\alpha
|:= \alpha_1+\cdots +\alpha_n$. If $|\alpha |=0$, i.e. $\alpha =0$,
then there is nothing to prove since $N_0= A^\d$. Suppose that $s:=
|\alpha |>0$ and the result is true for all $\alpha'$ with $|\alpha'
|<s$. Then, up to order, $\alpha_n>0$. Let $u\in N_\alpha$. We have
to show that $u\in N_\alpha'$. If $\d_n(u)=0$ then $u\in N_{\alpha -
e_n}$, and, by induction, $u\in N_{\alpha - e_n}\subseteq N_{\alpha
- e_n}'\subseteq N_\alpha'$ since $|\alpha - e_n|= |\alpha |
-1<|\alpha |$.

If $\d_n(u) \neq 0$ then $\d_n (u) \in N_{\alpha -e_n}= N_{\alpha -
e_n}'$ (by induction), and so $\d_n(u) = \sum_{0\leq \beta \leq
\alpha - e_n}\l_\beta x^{[\beta]}$ for some elements $\l_\beta \in
A^\d$, and not all of them  are equal to zero. Note that the element
$v:= \sum_{0\leq \beta \leq \alpha -e_n}\l_\beta x^{[\beta +e_n]}$
belongs to the set $N_\alpha'$, and $\d_n(u-v)=0$. Then $u-v\in
A^{\d_n}\cap N_\alpha = N_{\alpha'}$ where $\alpha':= (\alpha_1,
\ldots , \alpha_{n-1}, 0 )$. Since $|\alpha'|<|\alpha |$,
$N_{\alpha'}\subseteq N_{\alpha'}'$ and $v\in N_\alpha'$, we see
that $u= (u-v)+v\in N_\alpha'$, as required.

The second equality can be proved by similar arguments. $\Box $

$\noindent $

The following lemma gives all the $\d$-descents provided a single
$\d$-descent is known.

\begin{lemma}\label{a14Nov7}
Let $\d = (\d_1, \ldots , \d_n)$ be an $n$-tuple of commuting
derivations of a ring $A$, and $\{ x^{[\alpha]}, \alpha\in I:=
\prod_{i=1}^n [0,d_i]_{dis}\}$ be a $\d$-descent where $d_i\in \N
\cup \{ \infty \}$. Let $\{ x'^{[\alpha]}, \alpha \in I\}\subseteq
A$. Then the following statements are equivalent.
\begin{enumerate}
\item $\{ x'^{[\alpha]}, \alpha \in I\}$ is a $\d$-descent.
\item $x'^{[0]}:=1$ and, for each element $0\neq \alpha \in I$,
$x'^{[\alpha]}= x^{[\alpha]}+\sum_{0\neq \beta \leq \alpha}\l_\beta
x^{[\alpha -\beta]}$ where $\l_\g \in A^\d$, $0\neq \g \in I$.
\item $x'^{[0]}:=1$ and, for each element $0\neq \alpha \in I$,
$x'^{[\alpha]}= x^{[\alpha]}+\sum_{0\neq \beta \leq \alpha}
x^{[\alpha -\beta]}\mu_\beta$ where $\mu_\g \in A^\d$, $0\neq \g \in
I$.
\end{enumerate}
\end{lemma}

{\it Proof}. $(1\Leftarrow 2)$ This implication is obvious.

$(1\Rightarrow 2)$ Suppose that $\{ x'^{[\alpha]}, \alpha \in I\}$
is a $\d$-descent. Then $x'^{[0]}:=1$. Let $0\neq \alpha \in I$. By
Lemma \ref{a13Nov7} and the fact that $\d^\alpha (x^{[\alpha]})=1$,
we have $x'^{[\alpha]}=x^{[\alpha]}+\sum_{0\neq \beta \leq
\alpha}\l_\beta x^{[\alpha -\beta]}$ for some elements $\l_\beta \in
A^\d$ that depend on $\alpha$. For each $\g \in I$ such that $\g
\leq \alpha$,
$$ x'^{[\alpha - \g ]}= \d^\g (x'^{[\alpha ]}) = x^{[\alpha - \g
]}+\sum_{0\neq \beta \leq \alpha - \g } \l_\beta  x^{[\alpha - \beta
- \g ]},$$ and the implication $(1\Rightarrow 2)$  is obvious since
the equality above is true for all elements $\g , \alpha \in I$ such
that $\g \leq \alpha$.

$(1\Leftrightarrow 3)$ Repeat the arguments above making obvious
modifications.  $\Box $

\begin{lemma}\label{b14Nov7}
Let $\d = (\d_1, \ldots , \d_n)$ be an $n$-tuple of commuting
derivations of a ring $A$ such that $\d^\alpha (y^{[\alpha ]})=1$
for some elements $y^{[\alpha]}\in N(\d , A, \alpha )$, $\alpha
\in I:=\prod_{i=1}^n [0,d_i]_{dis}$ where $d_i\in \N \cup \{
\infty \}$. Then
\begin{enumerate}
\item there exists a unique $\d$-descent  $\{
x^{[\alpha]}, \alpha\in I \}$ such that,  for each nonzero element
$\alpha \in I$, $x^{[\alpha ]}= y^{[\alpha]}+\sum_{0\leq \beta
<\alpha}c_{\alpha \beta} y^{[\beta]}$ for some elements $c_{\alpha
\beta}\in A^\d$ such that $c_{\alpha , 0}=0$.
\item there exists a unique $\d$-descent  $\{
z^{[\alpha]}, \alpha\in I \}$ such that, for each nonzero element
 $\alpha \in I$, $z^{[\alpha ]}= y^{[\alpha]}+\sum_{0\leq \beta
<\alpha} y^{[\beta]}b_{\alpha \beta}$ for some elements $b_{\alpha
\beta}\in A^\d$ such that $b_{\alpha ,0}=0$.
\end{enumerate}
\end{lemma}

{\it Proof}. 1.  One can easily prove that $N_\alpha:= N(\d , A,
\alpha )=\oplus_{0\leq \beta \leq \alpha}A^\d y^{[\beta]}$,
$\alpha \in I$ (by modifying slightly the arguments of the proof
of Lemma \ref{a13Nov7}). Let us prove that there exists at least
one $\d$-descent, say $\{ x'^{[\alpha]}, \alpha \in I\}$. If the
set $I$ is finite then it contains the largest element with
respect to $\geq $, say $\beta \in I$. Then the multi-sequence $\{
x'^{[\alpha]}:=\d^{\beta - \alpha }(y^{[\beta ]}), \beta \in I\}$
is a $\d$-descent. If the set $I$ is infinite then we fix a
strictly ascending infinite sequence of elements of $I$,
 $\beta_1<\beta_2<\cdots $ such that $I=\cup_{i\geq 1} I_i$ where
 $I_i:= \{ \alpha \in I\, | \, \alpha \leq \beta_i\}$. For each
 $\alpha , \beta \in I$ such that $\alpha \leq \beta$, we have
 $$ \d^{\alpha}(y^{[\beta]})-y^{[\beta - \alpha]}\in
 \sum_{\g<\beta - \alpha }A^\d y^{[\g ]}.$$
Using this fact we can find elements $x'^{[\beta_i]}, i\geq 1$,
such that $x'^{[\beta_1]}:=y^{[\beta_1]}$ and
$\d^{\beta_i-\beta_{i-1}}(x'^{[\beta_i]})=x'^{[\beta_{i-1}]}$ for
all $i\geq 2$. Let $I_0:=\emptyset$. For each $\alpha \in I$,
there exists a unique index, say $i$, such that $\alpha \in
I_i\backslash I_{i-1}$. Let
$x'^{[\alpha]}:=\d^{\beta_i-\alpha}(x'^{[\beta_i]})$. Then it is
obvious that $\{x'^{[\alpha]}, \alpha \in I\}$ is the
$\d$-descent.

Let, for a moment, a multi-sequence  $ \{ x'^{[\alpha]}, \alpha \in
I\}$ be an arbitrary $\d$-descent in $A$. By Lemma \ref{a13Nov7},
for each $\alpha \in I$,  $x'^{[\alpha]}\in N_\alpha$, i.e.
$x'^{[\alpha]}=y^{[\alpha]}+\sum_{0\leq \beta
<\alpha}c_{\alpha\beta}'y^{[\beta]}$, for some elements
$c_{\alpha\beta}'\in A^\d$. Let $ \{ x^{[\alpha]}, \alpha \in I\}$
be another $\d$-descent, and so
$x^{[\alpha]}=y^{[\alpha]}+\sum_{0\leq \beta
<\alpha}c_{\alpha\beta}y^{[\beta]}$, $\alpha \in I$, for some
elements $c_{\alpha\beta}\in A^\d$. By Lemma \ref{a14Nov7},
$$x^{[\alpha]}=x'^{[\alpha]}+\sum_{0\neq \g \leq \alpha}\l_\g x'^{[\alpha - \g]}
\;\; {\rm for \; all}\;\;  \alpha \in I ,$$ for some elements $\l_\g
\in A^\d$.
 We have to prove that the defining conditions
$$c_{\alpha ,0}=0 \;\; {\rm for \; all}\;\;  0\neq \alpha \in I, $$ of the $\d$-descent from
Lemma \ref{b14Nov7}.(1) {\em uniquely} determine the elements
$\l_\g$. For each $i=1, \ldots , n$, the equality $c_{e_i,0}=0$
yields the equalities
$$y^{[e_i]}=x^{[e_i]}=x'^{[e_i]}+\l_{e_i}=y^{[e_i]}+c_{e_i,0}'+\l_{e_i},$$
hence $\l_{e_i}=-c_{e_i,0}'$. Let $s$ be a natural number such that
$s\geq 2$. Suppose that, using the equalities $c_{\alpha ,0}=0$,
$0\neq \alpha \in I$, $|\alpha |<s$,  we have already found unique
elements $\l_\g$ with $|\g |<s$. Take any element $\alpha \in I$
with $|\alpha | = s$. Then the element $\l_\alpha$ can be found
uniquely from the equality
$$x^{[\alpha]}=x'^{[\alpha]}+\sum_{0\neq \g <\alpha }\l_\g x'^{[\alpha - \g ]}+\l_\alpha .$$ For, we have to equate to zero the coefficient
$c_{\alpha ,0}$  of the elements $y^{[0]}:=1$ after we substitute
the sum for each $x'^{[\alpha - \g ]}$ above (via $y^{[\beta ]}$)
into the above equality:
$$
x^{[\alpha ]}=y^{[\alpha]}+\sum_{0\neq \beta
<\alpha}c_{\alpha\beta}y^{[\beta]}+(c_{\alpha,0}'+\sum_{0\neq \g
<\alpha} \l_\g c_{\alpha - \g , 0}'+\l_\alpha)y^{[0]},$$ that is
$\l_\alpha := -c_{\alpha ,0}'-\sum_{0\neq \g <\alpha} \l_\g
c_{\alpha - \g , 0}'$. Therefore,
 for this unique choice of the elements $\{ \l_\alpha \}$, we have $c_{\alpha ,0}=0$, $0\neq \alpha \in I$,
  for the $\d$-descent $\{ x^{[\alpha]}, \alpha \in I\}$. This proves the first statement of the lemma.

  2. Repeat the arguments of the first statement making obvious modifications.  $\Box$

$\noindent $

{\bf Structure of iterative sequence}. Structure of iterative
sequence of rank 1 is given by the following proposition.

\begin{proposition}\label{ittp}
{\rm (Structure of iterative sequence of rank 1, Proposition} 3.5,
\cite{Bav-simderharp}) Let $A$ be an $\Fp$-algebra and $\{ x^{[i]},
0\leq i <p^d \}$ be an iterative sequence. Then
\begin{enumerate}
\item for each $i=1, \ldots , p^d-1$, written $p$-adically as $i=
\sum_k i_kp^k$, $x^{[i]}= \prod_k \frac{x^{[p^k]i_k}}{i_k!}$. This
means that the iterative sequence is determined by the elements $\{
x^{[0]}, x^{[p^j]} \, | \,  j=0,1, \ldots ,  d-1\}$.  \item For each
$j=0,1,   \ldots, d-1$, $x^{[p^j]p}=0$ (hence $x^{[i]p}=0$ for all
$i=1, \ldots , p^d-1$, by statement 1). \item $x^{[0]}x^{[p^j]}=
x^{[p^j]}$, $j=0, 1, \ldots , n-1$, and $x^{[0]}x^{[0]}= x^{[0]}$.
\end{enumerate}
Conversely, given commuting elements $\{ x^{[0]}, x^{[p^j]}\, | \,
j=0,1, \ldots ,  d-1\}$, in $A$ that satisfy the conditions of
statements 2 and 3 above then the elements $\{ x^{[i]}, 0\leq i
<p^d\}$ defined as in statement 1 form an iterative sequence.
\end{proposition}

{\it Remark}. To make formulae more readable we often use the
notation $x^{[p^k]j}$ for $(x^{[p^k]})^j$.

Let $A$ be an $\Fp$-algebra. For each $i=1, \ldots , n$, let $\{
x_i^{[j]}, j\in [0,p^{d_i})_{dis}\}$ be an iterative sequence of
rank 1 in $A$ where $d_i\in \N \cup \{ \infty \}$. Suppose that
these sequences {\em commute} $(x_i^{[j]}x_k^{[l]}=
x_k^{[l]}x_i^{[j]})$ and have a common initial element, that is,
$x_1^{[0]}=x_2^{[0]}=\cdots = x_n^{[0]}$. It is easy to verify
that their product $\{ x^{[\alpha]}:=x_1^{[\alpha_1]}\cdots
x_n^{[\alpha_n]}, \alpha \in \prod_{i=1}^n[0, p^{d_i})_{dis}\}$ is
an iterative multi-sequence of rank $n$ and of exponent $(d_1,
\ldots , d_n)$.

$\noindent $

{\it Definition}. The iterative multi-sequence $\{
x^{[\alpha]}:=x_1^{[\alpha_1]}\cdots x_n^{[\alpha_n]}, \alpha \in
\prod_{i=1}^n[0, p^{d_i})_{dis}\}$ of rank $n$  and exponent
$(d_1, \ldots , d_n)$  is called the {\em product} of $n$
iterative sequences $\{ x_i^{[j]}, j\in [0,p^{d_i})_{dis}\}$ of
rank 1 and exponent $d_i$.

\begin{corollary}\label{c14Nov7}
{\rm (Structure of iterative multi-sequence of rank $n$)} Let $A$ be
an $\Fp$-algebra and $\{ x^{[\alpha]}, \alpha \in I \}$ be an
iterative multi-sequence of rank $n$ in $A$ where $I:=
\prod_{i=1}^n[0,p^{d_i})_{dis}$, $d_i\in \N \cup \{ \infty \}$. Then
the multi-sequence  $\{ x^{[\alpha]}, \alpha \in I \}$ is the
product of $n$ iterative sequences  $\{ x_i^{[j]}:=x^{[je_i]}, j\in
[0,p^{d_i})_{dis} \}$ of rank 1, and vice versa.
\end{corollary}

{\it Proof}. It is obvious that, for each $i=1, \ldots , n$, the
sequence $\{ x_i^{[j]}, j\in [0,p^{d_i})_{dis} \}$ is iterative of
rank 1. It is obvious that $x_1^{[0]}=x_2^{[0]}=\cdots = x_n^{[0]}$
and that $ x^{[\alpha]}:=x_1^{[\alpha_1]}\cdots x_n^{[\alpha_n]}$.
Now, the result follows. $\Box $

$\noindent $

{\bf Necessary and sufficient conditions for iterative
multi-sequence to be a $\d$-descent}. In case of sequence of rank 1
such conditions are given in the next proposition.

\begin{proposition}\label{d14Nov7}
{\rm (Corollary 3.6, \cite{Bav-simderharp})}  Let $A$ be an
$\Fp$-algebra, $\d$ be a derivation of $A$, and $\{ x^{[i]}, 0\leq
i<p^d\}$ be an iterative sequence in $A$ with $x^{[0]}=1$ where
$d\in \N \cup \{ \infty \}$. Then the iterative sequence $\{
x^{[i]}, 0\leq i<p^d\}$ is a $\d$-descent iff $\d (x^{[p^j]})=
x^{[p^j-1]}$, $0\leq j \leq d-1$.
\end{proposition}

\begin{corollary}\label{e14Nov7}
Let $\d = (\d_1, \ldots , \d_n)$ be an $n$-tuple of commuting
derivations of an $\Fp$-algebra  $A$, and $\{ x^{[\alpha]}, \alpha
\in I\}$ be an iterative multi-sequence of rank $n$ in $A$ with
$x^{[0]}=1$ where $I:=\prod_{i=1}^n [0,p^{d_i})_{dis}$, $d_i\in \N
\cup\{\infty\}$. Then the iterative multi-sequence $\{
x^{[\alpha]}, \alpha \in I\}$ is a $\d$-descent if and only if,
for each $i=1, \ldots , n$, the iterative sequence $\{
x_i^{[j]}:=x^{[je_i]},j\in [0,p^{d_i})_{dis}\}$ is a
$\d_i$-descent and $\{ x_i^{[j]},j\in [0,p^{d_i})_{dis}\}\subseteq
\cap_{k\neq i}A^{\d_k}$ if and only if
$\d_i(x_j^{[p^{k_j}]})=\d_{i,j}x_j^{[p^{k_j}-1]}$ for all $i,j=1,
\ldots , n$ and $k_j\in [0,d_j)_{dis}$ where $\d_{i,j}$ is the
Kronecker delta.
\end{corollary}

{\it Proof}. The first `if and only if' follows from Corollary
\ref{c14Nov7}. The second `if and only if' follows from Proposition
\ref{d14Nov7}. $\Box$

$\noindent $

Combining Proposition \ref{ittp}, Corollary \ref{c14Nov7} and
Corollary \ref{e14Nov7} we have necessary and sufficient conditions
for a multi-sequence to be an iterative $\d$-descent.

\begin{corollary}\label{f14Nov7}
Let $\d = (\d_1, \ldots , \d_n)$ be an $n$-tuple of commuting
derivations of an $\Fp$-algebra  $A$;  $\{ x_\nu^{[p^{k_\nu}]},
\nu = 1, \ldots , n; k_\nu=[0, d_\nu )_{dis}\}$ be commuting
elements of the algebra $A$ where $d_\nu \in \N \cup\{ \infty\}$;
$x_\nu^{[0]}:=1$ for all $\nu$. Let
 $x_\nu^{[i]}:=\prod_k  \frac{x_\nu^{[p^k]i_k}}{i_k!}$ for each
$i=\sum_k i_kp^k$, $0\leq i_k<p$, such that $0\leq i<p^{d_\nu}$. Let
$x^{[\alpha]}:=\prod_{\nu =1}^nx_\nu^{[\alpha_\nu]}$. Then the
 multi-sequence $ \{ x^{[\alpha]}, \alpha \in
I:=\prod_{i=1}^n[0,p^{d_\nu})_{dis}\}$  is an iterative $\d$-descent
iff $\d_\nu ( x_\mu^{[p^{k_\mu}]})=\d_{\nu , \mu}
x_\mu^{[p^{k_\mu}-1]}$ and $ x_\mu^{[p^{k_\mu}]p}=0$ for $\nu, \mu =
1,\ldots , n$ and $k_\mu \in [0,d_\mu)_{dis}$ (where $\d_{\nu ,
\mu}$ is the Kronecker delta).
\end{corollary}

Let $A$ be a commutative $\Fp$-algebra and $\d$ be a derivation of
$A$. Let $\ID (\d, d)$ be the set of all iterative $\d$-descents $
\{ x^{[i]}, 0\leq i <p^d\}$ of exponent $d$ in $A$ where $d\in \N
\cup\{ \infty\}$. Let $C(\d , d)$ be the set of all $d$-tuples
$(\l_0, \l_1, \ldots , \l_{d-1})$ such that $\l_j\in A^\d$ and
$\l_j^p=0$ for $0\leq j\leq d-1$. Note that if $A^\d$ is a {\em
reduced} ring then $C(\d , d) = \{ (0, \ldots , 0)\}$, i.e. $C(\d ,
d)$ contains a single element. In the case when $n=1$,
 by Lemma \ref{a13Nov7}  and Proposition \ref{ittp}, for each iterative $\d$-descent,
 say $ \{ x^{[i]}, 0\leq i
<p^d\}$, 
\begin{equation}\label{NdApn}
N(\d, A)_{p^d-1}= \bigoplus_{i=0}^{p^d-1}A^\d x^{[i]}\simeq A^\d [
x^{[1]}, x^{[p]}, \ldots , x^{[p^{d-1}]}]/(x^{[1]p}, x^{[p]p},
\ldots , x^{[p^{d-1}]p}).
\end{equation}
 So, $N(\d,
A)_{p^d-1}$ is the {\em subring} of $A$ that contains $A^\d$, and
the decomposition (\ref{NdApn}) holds for {\em all} iterative
$\d$-descents in $A$ of exponent $d$. In particular, all iterative
$\d$-descents in $A$ of exponent $d$ belong to the ring $N(\d,
A)_{p^d-1}$.
 Then
\begin{equation}\label{NNdA}
N(\d, A)_{p^d-1}=A^\d \oplus \gm, \;\; \gm :=(x^{[1]}, x^{[p]},
\ldots , x^{[p^{d-1}]}).
\end{equation}
If $\{ y^{[i]}, 0\leq i <p^d\}$ is an iterative $\d$-descent in $A$
then $\{ y^{[i]}, 0\leq i <p^d\}\subseteq N(\d, A)_{p^d-1}$.
Therefore, the following map is well-defined: 
\begin{equation}\label{rIDC}
r=r_d: \ID (\d , d)\ra C(\d , d), \;\; \{ y^{[i]}, 0\leq i
<p^d\}\mapsto (\l_0, \l_1, \ldots , \l_{d-1}),
\end{equation}
where $\l_j\equiv y^{[p^j]}\mod \gm$, $j=0, 1, \ldots , d-1$. Note
that the map $r$ depends on the choice of the iterative $\d$-descent
$\{ x^{[i]}, 0\leq i <p^d\}$ since the decomposition (\ref{NNdA})
does.

If $d=\infty$ then $N(\d , A)_{p^\infty -1}= N(\d, A)$, and so the
equalities (\ref{NdApn}) and (\ref{NNdA}) are as follows
\begin{eqnarray*}
 N(\d, A)&=& \bigoplus_{i\geq 0}A^\d x^{[i]}\simeq A^\d [
x^{[1]}, x^{[p]}, \ldots , x^{[p^i]}, \ldots ]/(x^{[1]p}, x^{[p]p},
\ldots , x^{[p^i]p}, \ldots )\\
N(\d, A)&=&A^\d \oplus \gm, \;\; \gm :=(x^{[1]}, x^{[p]}, \ldots ,
x^{[p^i]}, \ldots ).
\end{eqnarray*}
Let $\ID (\d ) := \ID (\d , \infty )$. Then the map (\ref{rIDC})
takes the form 
\begin{equation}\label{8inf}
r=r_\infty : \ID (\d )\ra C(\d , \infty), \;\; \{ y^{[i]},  i\in \N
\}\mapsto (\l_0, \l_1, \ldots , \l_j, \ldots ),
\end{equation}
where $\l_j\equiv y^{[p^j]}\mod \gm$,  $j\in \N$.

The following theorem is a key (and difficult) result of the paper
\cite{Bav-simderharp}.

\begin{theorem}\label{f18Jan06}
{\rm (Existence and uniqueness of an iterative $\d$-descent when
$n=1$, Theorem 3.8, \cite{Bav-simderharp})} Let $A$ be a commutative
algebra over a field $K$ of characteristic $p>0$ and $\d $ be a
$K$-derivation of the algebra $A$ such that there exists a finite
sequence of elements $y_0, y_1, \ldots , y_{d-1}$ of $A$ such that
$y_k^p=0$ and $\d^{p^k}(y_k)=1$ for all $0\leq k\leq d-1$. Then
\begin{enumerate}
\item (Existence) The following sequence $\{ x^{[i]}, 0\leq i
<p^d\}$ is an iterative $\d$-descent where $x^{[0]}:=1$,
$x^{[1]}:=y_0$,  and, for  $i\geq 2$ written $p$-adically as
$i=\sum_{k=0}^ti_kp^k$ ($0\leq i_k\leq p-1$) the element $x^{[i]}$
is defined as
$x^{[i]}:=\prod_{k=0}^t\frac{(x^{[p^k]})^{i_k}}{i_k!}$, where
$$x^{[p]}:= (-1)^{p-1}\phi_0(y_1), \;\; \phi_0 (z):=
\sum_{j=0}^{p-1}(-1)^j\frac{(x^{[1]})^j}{j!}\d^j(z),$$ and then
recursively, for each $k$ such that  $1\leq k \leq n-2$, the element
$x^{[p^{k+1}]}$ is defined by the rule
$$x^{[p^{k+1}]}:= (-1)^{p-1}\d^{p^k-1} (\prod_{l=0}^{k-1}
\frac{(x^{[p^l]})^{p-1}}{(p-1)!}\cdot \phi_k(y_{k+1})), \;\;
\phi_k (z):=
\sum_{j=0}^{p-1}(-1)^j\frac{(x^{[p^k]})^j}{j!}\d^{p^kj}(z).$$
\item (Almost uniqueness) Let $\{ x^{[i]}, 0\leq i <p^d\}$ be an
arbitrary  iterative $\d$-descent (not necessarily as in statement
1, and $d$ here is not necessarily as in statement 1 either). Then
the map (\ref{rIDC}) is a bijection. \item (Uniqueness). If, in
addition, the ring $A^\d$ is reduced then $\{ x^{[i]}, 0\leq i
<p^d\}$ from statement 1  is the only iterative $\d$-descent.
\end{enumerate}
\end{theorem}

{\it Remark}. Note that Theorem \ref{f18Jan06} holds also for
$d=\infty$ since each infinite iterative $\d$-descent $\{ x^{[i]},
i\in \N\}$ is a union of finite iterative $\d$-descents $\{ x^{[i]},
i\in [0,p^d)_{dis}\}$, $d\in \N$.

$\noindent $

Let $A$ be an $\Fp$-algebra and $\d = (\d_1, \ldots , \d_n)$ be an
$n$-tuple of commuting derivations of the algebra $A$. Let $\ID (\d
, d )$ be the set of all the iterative $\d$-descents $\xx := \{
x^{[\alpha]}\}$ of exponent $d:=(d_1, \ldots , d_n)$ where  $d_i\in
\N \cup \{ \infty \}$. Let $C(\d , d ):=C'( d_1 )\times \cdots
\times C'(d_n )$ where $C'(d_i)$ is the set of all $d_i$-tuples
$(\l_0, \l_1, \ldots , \l_{d_i-1})$ such that $\l_j\in A^\d$ and
$\l_j^p=0$ for all $j\in [0,d_i)_{dis}$. By Corollary \ref{c14Nov7},
each $\xx$ is the product $\prod_{i=1}^n\xx_i$ of the iterative
$\d_i$-descents $\xx_i:=\{ x_i^{[j]}:=x^{[je_i]}, j\in [0,
p^{d_i})_{dis}\}$. For each $i=1, \ldots , n$, let $r_{d_i}:\ID
(\d_i, d_i )\ra C(\d , d_i )$ be the map (\ref{rIDC}) for the
derivation $\d_i$. Consider the map 
\begin{equation}\label{a8inf}
r_d: \ID (\d , d)\ra C(\d , d), \;\; \xx \mapsto (r_{d_1}(\xx_1),
\ldots , r_{d_n}(\xx_n)).
\end{equation}

\begin{theorem}\label{15Nov7}
{\rm (Existence and uniqueness of an iterative $\d$-descent)} Let
$A$ be a commutative $\Fp$-algebra and $\d =(\d_1, \ldots , \d_n)$
be
 an $n$-tuple of commuting derivations of the algebra $A$. Suppose
 that there exists a set of elements of the algebra $A$, $\{ y_{ik},
 i=1, \ldots , n; k\in [0,d_i)_{dis}\}$, where $d_i\in \N \cup\{\infty\}$,  such that $y_{ik}^p=0$ and
 $\d_i^{p^k}(y_{ik})=1$ for all $i=1, \ldots , n$ and $k\in
 [0,d_i)_{dis}$, and $\d_j(y_{ik})=0$ for all $i\neq j$. Note that, for each $i=1,
 \ldots , n$, the  elements $\{y_{ik}, k\in [0,d_i)_{dis}\}$
 satisfy the assumptions of
 Theorem \ref{f18Jan06}.(1), let $\xx_i:=\{ x_i^{[j]}, j\in
 [0,p^{d_i})_{dis}\}$ be the corresponding $\d_i$-descent of rank 1 and  exponent $d_i$ as in Theorem
 \ref{f18Jan06}.(1). Then
\begin{enumerate}
\item (Existence) the product $\xx :=\prod_{i=1}^n \xx_i=\{
x^{[\alpha]}:=x_1^{[\alpha_1]}\cdots x_n^{[\alpha_n]}, \alpha \in
\N^n\}$ of the iterative $\d_i$-descents $\xx_i$ is an iterative
$\d$-descent.
\item (Almost uniqueness) The map (\ref{a8inf}) is a bijection.
\item (Uniqueness). If, in addition, the ring $A^\d$ is
reduced then the iterative $\d$-descent $\xx$ from statement 1 is
the only iterative $\d$-descent in $A$.
\end{enumerate}
\end{theorem}

{\it Proof}. 1. Statement 1 follows from Corollary \ref{d14Nov7}.

2. Statement 2 follows from statement 1 and Theorem
\ref{f18Jan06}.(2).

3.  Statement 3 follows from statement 2 since the set $C(\d , d)$
contains a single element (since the ring $A^\d$ is reduced). $\Box
$

$\noindent $


\section{Rigidity of the group $\Aut_K(\CDPn )$}\label{RIGID}

In this section, the rigidity of the group $\Aut_K(\CDPn )$ is
proved (Theorem \ref{u2Jul7}).

$\noindent $

{\bf The ring $\CD (P_n)$  of  differential operators}.  The ring
$\CD (P_n)$ of differential operators on a polynomial algebra $P_n:=
K[x_1, \ldots , x_n]$ is a $K$-algebra generated by the elements
$x_1, \ldots , x_n$ and {\em commuting} higher derivations $\derik
:=\frac{\der_i^k}{k!}$, $i=1, \ldots , n$; $k\geq 1$,  that satisfy
the following defining relations: 
\begin{equation}\label{DPndef}
[x_i,x_j]=0, \;\; [\derik , \derjl ]=0,\;\;\; \derik \deril
={k+l\choose k}\der_i^{[k+l]}, \;\;\; [\derik ,
x_j]=\d_{ij}\der_i^{[k-1]},
\end{equation}
 for all
$i,j=1, \ldots , n$;  $k,l\geq 1$, where $\d_{ij}$ is the
Kronecker delta,  $\der_i^{[0]}:=1$, $\der_i^{[-1]}:=0$, and
$\der_i^{[1]}=\der_i=\frac{\der}{\der x_i}\in \Der_K(P_n)$,
$i=1,\ldots , n$. The action of the higher derivation $\derik$ on
 the polynomial algebra $P_n=K\t_\Z \Z [x_1, \ldots , x_n]$ should
 be understood as the action of the element $1\t_\Z
 \frac{\der_i^k}{k!}$.

The algebra $\CD (P_n)$ is a {\em simple} algebra. Note  that the
algebra $\CD (P_n)$ is  not finitely generated and neither left
nor right Noetherian, it  does not satisfy finitely many defining
relations.
$$\CD (P_n)= \bigoplus_{\alpha , \beta \in \Nn }Kx^\alpha \derbb
=\bigoplus_{\alpha , \beta \in \Nn }K \derbb x^\alpha
  , \;\; {\rm where} \;\; x^\alpha := x_1^{\alpha_1}\cdots
x_n^{\alpha_n}, \; \derbb := \der_1^{[\beta_1]}\cdots
\der_n^{[\beta_n]}.$$ For each $i=1, \ldots , n$, let $\CD
(P_1)(i):= \CD (K[x_i])$. Then $$\CDPn = \t_{i=1}^n \CD
(P_1)(i)\simeq \CD(P_1)^{\t n}.$$ For each $i=1, \ldots , n$ and
$j\in \N$ written $p$-adically as $j= \sum_k j_kp^k$, $0\leq j_k<p$,
\begin{equation}\label{Fdij}
\derij=\prod_k\derijkpk = \prod_k \frac{\deripkjk}{j_k!}, \;\;
\derijkpk =\frac{\deripkjk}{j_k!},
\end{equation}
where $\deripkjk := (\deripk )^{j_k}$. For $\alpha , \beta \in \Nn$,
\begin{equation}\label{daaxb}
\derba (x^\beta ) = {\beta \choose \alpha } x^{\beta - \alpha}, \;\;
{\beta \choose \alpha }:= \prod_i {\beta_i \choose \alpha_i},
\end{equation}
where, in the formula above, $x_i^t:=0$ for all negative integers
$t$ and all $i$.

For $\alpha , \beta \in \Nn$, 
\begin{equation}\label{dadb}
\derba \derbb = {\alpha +\beta \choose \beta}\der^{[\alpha +\beta
]}.
\end{equation}

The algebra $\CDPn$ has two natural filtrations: the {\em
canonical filtration}  $F=\{ F_i\}_{i\geq 0}$ where $F_i:=
\oplus_{|\alpha |+|\beta |\leq i}Kx^\alpha \der^{[\beta]}$ and
 $|\alpha |:=|\alpha_1|+\cdots + |\alpha_n|$, and the {\em order filtration} $\{ \CDPn_i\}_{i\geq 0}$
where $\CDPn_i:= \oplus_{|\beta |\leq i} P_n\der^{[\beta ]}$. The
first one is a finite dimensional filtration but the second one is
not.
 For both filtrations the associated graded algebras are
 isomorphic to the tensor product $P_n\t \D_n$ of the polynomial
 algebra $P_n$ and the algebra $\D_n=\oplus_{\alpha \in \N^n}
 K\der^{[\alpha}]$ of scalar differential operators. In
 particular, these algebras are commutative but not finitely
 generated, not Noetherian, and  not domains.

$\noindent $

{\bf Defining relations for $\CD (P_n)$}. As an abstract
$K$-algebra, the ring $\CD (P_n)$ of differential operators on $P_n$
is generated by the elements $x_i$, $\deripk$, $i=1, \ldots , n$,
$k\in \N$, that satisfy the defining relations: 
\begin{equation}\label{RelsDPn}
[x_i,x_j]=0, \;\; [\deripk , \derjpl ]=0, \; (\deripk )^p=0, \;
[\deripk , x_j]=
\d_{ij}\prod_{l=0}^{k-1}\frac{\der_i^{[p^l](p-1)}}{(p-1)!}
\end{equation}
for all $1\leq i,j\leq n$  and  $k,l\in \N$.

In  characteristic $p>0$, the next theorem proves existence and
uniqueness of iterative $(-\ad (x_1), \ldots ,$ $  -\ad
(x_n))$-descent in $\CDPn$. In characteristic zero, this result is
not true, there are many iterative $(-\ad (x_1), \ldots , -\ad
(x_n))$-descents in $\CDPn$: $\{ y^{[\alpha ] }:=
\prod_{i=1}^n\frac{(\der_i+a_i)^{\alpha_i}}{\alpha_i!}, \alpha \in
\Nn \}$ is the iterative $(-\ad (x_1), \ldots , -\ad
(x_n))$-descent in $\CDPn$ where $a_i\in K[x_i]$, $i=1, \ldots ,
n$.

Since the ring of invariants $\cap_{i=1}^n \D_n^{\ad (x_i)}=K$ is
reduced, the multi-sequence $\{ \der^{[\alpha]}, \alpha \in \Nn
\}$ is the only iterative $(-\ad (x_1), \ldots , -\ad
(x_n))$-descent in the algebra $\D_n$, by Theorem \ref{15Nov7}.
The next theorem proves the same but for the ring $\CDPn$. The
ring $\CDPn$ is non-commutative, so we cannot apply Theorem
\ref{15Nov7} directly.

\begin{theorem}\label{u2Jul7}
Let $K$ be a field of characteristic $p>0$. Then $\{
\der^{[\alpha]}, \alpha \in \Nn \}$ is the only iterative $(-\ad
(x_1), \ldots , -\ad (x_n))$-descent in $\CDPn$.
\end{theorem}

{\it Proof}.   Let $\d_1:=   - \ad (x_1), \ldots , \d_n:=-\ad
(x_n)$ and $\d := (\d_1,\ldots , \d_n)$. Let $\{ y^{[\alpha]},
\alpha\in \Nn\}$ be an iterative $\d$-descent in $\CDPn$. We have
to show that $y^{[\alpha]}= \der^{[\alpha]}$ for all $\alpha$ or,
equivalently, that $y_i^{[p^k]}= \der_i^{[p^k]}$ for all $i=1,
\ldots , n$ and $k\in \N$.  To prove this fact we use induction on
$k$.

Let $k=0$ and $y_i:=  y_i^{[p^0]}$. Let   us prove that $y_i=
\der_i$ for all $i$.  Note that
$$ y_i\in \cap_{j\neq i} \ker ( \d_j) \cap  \ker(\d_i^2) =
P_n\oplus P_n\der_i.$$ Since $\d_i(y_i) = 1$, we see that $y_i=
\der_i+\l_i$ for  some polynomial $\l_i\in  P_n$. If $\l_i=0$, we
are done. Suppose that $\l_i\neq 0$, we seek a contradiction.
Then, using the defining relations (\ref{RelsDPn})  for $\CDPn$,
we obtain that
$$ 0=y_i^p=(\der_i+\l_i)^p= \der_i^p+\sum_{j=1}^{p-1} a_j\derij
+a_0+\l_i^p$$ for some polynomials $a_k\in P_n$ such that
$\deg_{x_i}(a_0)<\deg_{x_i}(\l_i^p)$ where $\deg_{x_i}(q)$ is the
$x_i$-degree of a polynomial $q\in P_n$. Since $\der_i^p=0$, the
above equality yields $a_1=\cdots =a_{p-1}=0$ and $a_0+\l_i^p=0$.
The last equality is impossible since
$\deg_{x_i}(a_0)<\deg_{x_i}(\l_i^p)$.

Suppose that $k\geq 1$ and that $y_i^{[p^s]}=\der_i^{[p^s]}$ for
all $s=0, \ldots , k-1$ and  $i=1, \ldots , n$. To finish the
proof by induction it remains to show  that $y_i^{[p^k]}=
\der_i^{[p^k]}$ for all $i$. Note that

\begin{eqnarray*}
 y_i^{[p^k]}&\in & \cap_{j\neq i} \ker (\d_j) \cap \ker (\d_i^{[p^{k+1}]})=
 P_n\oplus K[\der_i, \der_i^{[p]}, \ldots , \deripk ]\;\;\;\;\; ({\rm Lemma}\;\; \ref{a28Jun7}.(6)), \\
 y_i^{[p^k]}&\in & \cap_{j=1}^n \ker \, \ad (\der_j^{[p^{k-1}]}))=
  K[x_1^{p^k}, \ldots , x_n^{p^k} ]\t \D_n \;\;\;\;\;\;\;\;\;\;\;\;\;\;\;\;\;\;\, ({\rm Lemma}\;\; \ref{a28Jun7}.(2)),
\end{eqnarray*}
and so $y_i^{[p^k]}\in K[x_1^{p^k}, \ldots , x_n^{p^k} ]\t
K[\der_i, \der_i^{[p]}, \ldots , \deripk ]$. Since
$\d_i^{[p^k]}(y_i^{[p^k]})=1$,
$$ y_i^{[p^k]}= \deripk + \sum_{j=0}^{p^k-1}\mu_j\derij$$
 for some polynomials $\mu_j\in K[x_1^{p^k}, \ldots , x_n^{p^k}
 ]$. On the one hand, $[y_i^{[p^k]}, x_i]=y_i^{[p^k-1]}=
 \der_i^{[p^k-1]}$. On the other hand,
 $$[y_i^{[p^k]},
 x_i]=\der_i^{[p^k-1]}+\sum_{j=1}^{p^k-1}\mu_j\der_i^{[j-1]}.$$
 Hence, $\mu_1=\cdots = \mu_{p^k-1}=0$, i.e. $y_i^{[p^k]}=
 \der_i^{[p^k]}+\mu_0$. If $\mu_0=0$, we are done. Suppose that
 $\mu_0\neq 0$, we seek a contradiction. Then
$$0=(y_i^{[p^k]})^p=(\deripk +\mu_0)^p=(\deripk )^p
+\sum_{j=1}^{p^k-1}b_j\derij+b_0+\mu_0^p$$ for some polynomials
$b_s\in K[x_1^{p^k}, \ldots x_n^{p^k}]$ such that
$\deg_{x_i^{p^k}}(b_0)<\deg_{x_i^{p^k}}(\mu_0^p)$. Since
$(\deripk)^p=0$,  the above equality yields $b_1=\cdots
=b_{p^k-1}=0$ and $b_0+\mu_0^p=0$. The last equality is impossible
since $\deg_{x_i^{p^k}}(b_0)<\deg_{x_i^{p^k}}(\mu_0^p)$. This
proves that $y_i^{[p^k]}= \deripk$ for all $i$ and $k\in \Nn$ The
proof of the theorem is complete. $\Box $

$\noindent $

Each automorphism $\s \in \Aut_K(P_n)$ can be naturally extended
(by change of variables) to a $K$-automorphism, say $\s$, of the
ring $\CDPn$ of differential operators on the polynomial algebra
$P_n$ by the rule: 
\begin{equation}\label{hs=sas}
\s(a) := \s a \s^{-1}, \;\; a\in \CDPn .
\end{equation}
Then the group $\Aut_K(P_n)$ can be seen as a subgroup of $\Aut_K
(\CDPn )$ via (\ref{hs=sas}).  In characteristic zero, there are
many other extensions of the automorphism $\s$. For example, the
automorphisms of $\CD (K[x])$:
$$\s_f : x\mapsto x, \;\; \der \mapsto \der +f, \;\; f\in K[x],$$
are extensions of the identity automorphism of $K[x]$. This is not
the case in prime characteristic as Theorem \ref{28Jun7} shows.

The set $\{ x_i, \deripk \, | \, 1\leq i \leq n, k\in \N \}$ is a
 set of $K$-algebra generators for $\CDPn$. Each automorphism $\s
$ of $\CDPn$ is uniquely determined by its action on this {\em
infinite} set. Theorem \ref{28Jun7} (which is not true in
characteristic zero) says that $\s$ is uniquely determined by the
elements $\s (x_1), \ldots , \s (x_n)$.

$\noindent $

{\bf Proof of Theorem \ref{28Jun7}}. $(\Rightarrow )$ This
implication is trivial.

$(\Leftarrow )$ By considering the automorphism $\s\tau^{-1}$ of
the algebra $\CDPn$, it suffices to show that if  an automorphism,
say $\s$, of the algebra $\CDPn$ satisfies the conditions that $\s
(x_1)=x_1, \ldots , \s (x_n) = x_n$, then $\s$ is the identity
map. Using the defining relations (\ref{DPndef}), we see that $\{
\s (\der^{[\alpha]}), \alpha \in \N^n\}$ is the iterative $(-\ad
(x_1), \ldots , -\ad (x_n))$-descent. By Theorem \ref{u2Jul7}, $\s
(\der^{[\alpha]})= \der^{[\alpha]}$ for all elements $\alpha \in
\N^n$, i.e. the automorphism $\s$ is the identity map, as
required. $\Box $


\begin{corollary}\label{a16Mar8}
$\Aut_K(P_n)= \{ \tau \in \Aut_K(\CDPn )\, | \, \tau (P_n) =
P_n\}$ and, for each automorphism $\s \in \Aut_K(P_n)$, the
equation (\ref{hs=sas}) gives the only extension of the
automorphisms $\s$ to an automorphism of the $K$-algebra $\CDPn$.
\end{corollary}


\section{Extensions of the Frobenius to the ring $\CDPn$}\label{ETFRB}

In this section, Theorem \ref{18Mar8} is proved and the concept of
a Frobenius on $\CDPn$ is introduced.

$\noindent $

{\bf The canonical extension of the Frobenius $F$ to $\CD (P_n)$}.
Let $A$ be an algebra over a field $K$ of characteristic $p>0$.
For each $x\in A$, 
\begin{equation}\label{adxp}
({\rm ad} \, x)^p= {\rm ad }  (x^p)
\end{equation}
since, for any $ a\in A$, ${\rm ad }  (x^p)(a)= [ x^p,
a]=\sum_{i=1}^p{p\choose i} ({\rm ad}\, x)^i (a)x^{p-i}=({\rm
ad}\, x)^p(a)$.

Using the defining relations (\ref{RelsDPn}) for the algebra
$\CDPn $ and the equality (\ref{adxp}), we see that  the {\em
Frobenius} $\Fp$-algebra monomorphism
$$F:P_n\ra P_n, \;\; a\mapsto a^p,$$ can be
lifted to the  $\Fp$-algebra monomorphism of the ring $\CDPn$  by
the rule
\begin{equation}\label{Flift}
F=F_x:\CD (P_n)\ra \CD (P_n), \;\; \derba \mapsto \der^{[p\alpha
]}, \;\ \alpha \in \Nn.
\end{equation}
In more detail, it suffices to check that
$$ [\der_i^{[p^{k+1}]}, x_i^p]=
\prod_{l=1}^{k}\frac{\der_i^{[p^l](p-1)}}{(p-1)!}, \;\; 1\leq i
\leq n, \;\; k\geq 0.$$ By (\ref{adxp}) and (\ref{daaxb}),
\begin{eqnarray*}
 [\der_i^{[p^{k+1}]}, x_i^p]&=& (-\ad \, x_i)^p (\der_i^{[p^{k+1}]})=
 (-\ad \, x_i)^{p-1} (\prod_{l=0}^{k}\frac{\der_i^{[p^l](p-1)}}{(p-1)!}) \\
 &=& (-\ad \, x_i)^{p-1} (\frac{\der_i^{p-1}}{(p-1)!})\cdot
 (\prod_{l=1}^{k}\frac{\der_i^{[p^l](p-1)}}{(p-1)!})=\prod_{l=1}^{k}\frac{\der_i^{[p^l](p-1)}}{(p-1)!}.
\end{eqnarray*}

{\it Definition}. The $\Fp$-algebra monomorphism $F=F_x$ in
(\ref{Flift}) is called the {\bf canonical Frobenius} with respect
to the choice $x=(x_1, \ldots , x_n)$ of  free generators for the
polynomial algebra $P_n$.

$\noindent $

Note that the Frobenius $F$ in (\ref{Flift}) is {\em not} the
`obvious' extension of the Frobenius (which makes no sense as it
can be easily seen).

$\noindent $

{\it Definition}. A ring homomorphism $F': \CDPn \ra \CDPn$ is
called a {\bf Frobenius} (homomorphism) if the following
conditions hold:

1. there exists an $F'$-invariant polynomial subalgebra $P_n'=
K[x_1', \ldots , x_n']$ of the ring $\CDPn$ such that the
restriction $F'|_{P_n'}$ is the Frobenius on $P_n'$: $a\mapsto
a^p$;

2. the inner derivations $\ad (x_1'), \ldots , \ad (x_n')$ of the
ring $\CDPn$ are locally nilpotent derivations with $\cap_{i=1}^n
\CDPn^{\ad (x_i')}= P_n'$;

3. there exist a commutative $F'$-invariant $K$-subalgebra $\D'$
of the ring $\CDPn$ such that $[x_i', \D']\subseteq \D'$ for all
$i$ and elements $y_1, \ldots , y_n\in \D'$ such that
$y_1^p=\cdots = y_n^p=0$, $[y_i, x_j']=\d_{ij}$ (the Kronecker
delta), and $[F'^k(y_i) , x_j']=0$ for all $i\neq j$ and $k\geq
1$.

$\noindent $

{\it Example 1}. The canonical Frobenius $F_x$ is a Frobenius
where $x_i'=x_i$, $\der_i'=\der_i$ and $\D'= \D_n$, the algebra of
scalar differential operators.

$\noindent $

{\it Example 2}. For each $K$-algebra automorphism $\s \in
\Aut_K(\CDPn)$, $F':= \s F_x\s^{-1}$ is a Frobenius on $\CDPn$
where $x_i'= \s (x_i)$, $\der_i= \s (\der_i)$ and $\D' = \s
(\D_n)$.

$\noindent $

 Let $\Frob (\CD (P_n))$ be the set of all the Frobenius
homomorphisms on $\CDPn$. The polynomial algebra $P_n'$ above is
called a {\em Frobenius polynomial subalgebra}  of the ring
$\CDPn$. Let $\FPolCDPn$ be the set of all the  Frobenius
polynomial subalgebras of the ring $\CDPn$. We will see later that
each Frobenius subalgebra is a maximal commutative subalgebra of
the algebra $\CDPn$ (Corollary \ref{a21Mar8}.(1)), and so one
Frobenius polynomial subalgebra cannot be a proper subalgebra of
another. We will see later that the algebra $\D'$ is unique for
each Frobenius $F'$ on $\CDPn$  and it is also a maximal
commutative subalgebra of $\CDPn $ (Corollary
\ref{a21Mar8}.(1,3)). Let $\Frob (\CDPn , P_n')$ be the set of all
the Frobenius homomorphisms $F'$ with the Frobenius polynomial
subalgebra $P_n'$. Then 
\begin{equation}\label{FruP}
\Frob (\CDPn )= \bigcup_{P_n'\in \FPolCDPn } \Frob (\CDPn , P_n').
\end{equation}
This  is  the disjoint union since each Frobenius polynomial
subalgebra is a maximal commutative subalgebra. The group
$\Aut_K(\CDPn )$ acts naturally on the sets $\Frob (\CDPn )$ and
$\FPolCDPn$: $(\s , F')\mapsto \s F'\s^{-1}$ and $(\s ,
P_n')\mapsto \s (P_n')$. Clearly, $\s \Frob (\CDPn , P_n')\s^{-1}
= \Frob (\CDPn , \s (P_n'))$.

$\noindent $

{\it Definition}. $\Frob (\CDPn, P_n, \D_n):= \{ F'\in \Frob
(\CDPn , P_n)\, | \, F'(\D_n) \subseteq \D_n\}$.

$\noindent $

{\bf  Proof of Theorem \ref{18Mar8}}. Let $F'\in \FrobCDPn$ and
let $x_i'$, $y_i$ and $\D'$ be as in the Definition of Frobenius
homomorphism on $\CDPn$ above. Let $y_{ik}:= F'^k(y_i)$, $k\geq
0$, and $\d_i:= - \ad (x_i')$. Then the conditions of Theorem
\ref{15Nov7}.(3) hold for the algebra $A:= \D'$, $\d =(\d_1,
\ldots , \d_n)$ and the elements $\{ y_{ik}\}$:
\begin{eqnarray*}
 y_{ik}^p&=&F'^k(y_i)^p= F'^k(y_i^p)=0;  \\
 \d_i^{p^k}(y_{ik})&=&-\ad (x_i')^{p^k} (y_{ik})= -\ad
 (x_i'^{p^k})(y_{ik})= [ y_{ik}, x_i'^{p^k}]=[F'^k(y_i),
 F'^k(x_i')]\\
 &=&F'^k([y_i, x_i'])= F'^k(1) = 1;\\
 \d_j(y_{ik})&=& [ F'^k(y_i), x_j']=0, \;\; i\neq j, \;\; k\geq 0;
\end{eqnarray*}
the algebra $\D'^\d = \D'\cap \CDPn^\d = \D'\cap P_n'$ is reduced.
By Theorem \ref{15Nov7}.(3), in $\D'$ there exists a unique
iterative $\d$-descent, say $\{ \der'^{[\alpha]}\}$. Since the
derivations $\d_i$ of the algebra $\CDPn$ are locally nilpotent
and $\CDPn^\d = P_n'$, by using Lemma \ref{a13Nov7} we have
$$ \CDPn = N(\d , \CDPn ) = \bigoplus_{\alpha \in \N^n} \CDPn^\d
\der'^{[\alpha]}= \bigoplus_{\alpha\in
\N^n}P_n'\der'^{[\alpha]}\simeq \CD (P_n') \simeq \CDPn . $$
Consider the $K$-algebra automorphism $\s$ of $\CDPn$ given by the
rule: $x_i\mapsto x_i'$, $\der^{[\alpha]}\mapsto
\der'^{[\alpha]}$. Then $\s (P_n) = P_n'$, and so $\FPolCDPn =
\Aut_K(\CDPn ) \cdot P_n$. The algebra $\D_n=\oplus_{\alpha \in
\N^n}K\der^{[\alpha]}$ of scalar differential operators is a {\em
maximal} commutative subalgebra of $\CDPn$, hence $\oplus_{\alpha
\in \N^n} K\der'^{[\alpha]}$ is the  maximal commutative
subalgebra of $\CDPn$ which is contained in the commutative
subalgebra $\D'$. Therefore, $\D' = \oplus_{\alpha \in \N^n}
K\der'^{[\alpha]}$. This means that $\s (\D_n)=\D'$, and so
$\s^{-1} F'\s \in \Frob (\CDPn , P_n, \D_n)$ and the algebra $\D'$
in the Definition of a Frobenius $F'$ of $\CDPn$ is unique, and
it is a maximal commutative subalgebra of $\CDPn$. This proves
statement 1 and Corollary \ref{a21Mar8}.(1,2)).

The stabilizer of the polynomial algebra $P_n$ under the action of
the group $\Aut_K(\CDPn )$ on the set $\FPolCDPn$ is equal to the
set $\{ \s\in \Aut_K(\CDPn ) \, | \, \s (P_n) = P_n\}=
\Aut_K(P_n)$, Corollary \ref{a16Mar8}. Now, $\Aut_K(\CDPn ) \cdot
P_n \simeq \Aut_K(\CDPn ) / \Aut_K(P_n)$. This completes the proof
of statement 3.

To prove statement 2, let $F'\in \Frob (\CDPn , P_n)$. By
statement 1, $F' = \s F''\s^{-1}$ for some $F''\in \Frob (\CDPn ,
P_n, \D_n)$ and $\s \in \Aut_K(\CDPn )$ such that $\s (P_n) =
P_n$. Then $\s \in \Aut_K(P_n)$, by Corollary \ref{a16Mar8}, as
required. $\Box $

\begin{corollary}\label{a21Mar8}

\begin{enumerate}
\item For each Frobenius $F'\in \Frob (\CDPn )$, the corresponding
algebras $P_n'$ and $\D'$ are unique, they are maximal commutative
subalgebras of the algebra $\CDPn$. \item Let $F', F''\in
\Frob(\CDPn )$ and $F'' = \s F'\s^{-1}$ for some automorphism $\s
\in \Aut_K(\CDPn )$. Then $P_n''= \s (P_n')$ and $\D'' = \s (\D')$
 where $P_n''$ and $\D''$ are the corresponding algebras for the Frobenius $F''$.
 \item One Frobenius polynomial algebra cannot be a proper
 subalgebra of another Frobenius polynomial algebra. \item The
 union (\ref{FruP}) is a disjoint union.
\end{enumerate}
\end{corollary}


\section{The sets $\FrobsCDPnDn $, $s\geq 1$}

$\noindent $

{\it Definition}. For each natural number $s\geq 1$, let $
\FrobsCDPnPn$ be the set  of all  ring  endomorphisms $F'$ of
$\CDPn$ such that $F'(P_n)  \subseteq   P_n$ and $F'|_{P_n}= F^s$,
i.e. $F'(a) = a^{p^s}$ for all $a\in P_n$.

$\noindent $

{\it Definition}. For each natural number $s\geq 1$, let
$$
\FrobsCDPnDn:= \{  G\in \FrobsCDPnPn\, | \, G(\D_n)\subseteq
\D_n\}.$$

$\noindent $

This set contains  precisely all   the extensions $G$ of $F^s$ on
the polynomial algebra $P_n$ that respect the subalgebra $\D_n $
of scalar differential operators of $\CDPn$.  In this section, the
sets $\FrobsCDPnDn$ are found explicitly (Theorem \ref{M5Jul7} and
Theorem \ref{cM5Jul7}).


$\noindent $

Let $\d := (\d_1, \ldots , \d_n)$ where $\d_i:= -\ad (x_i)$ is the
inner derivation of the ring $\CDPn$. For each natural number $s\geq
1$, $\d_i^{p^s}=-\ad (x_i)^{p^s}= -\ad (x_i^{p^s})$ is the inner
derivation of the ring $\CDPn$. For each natural number $s\geq 1$,
$\d^{p^s}:=(\d_1^{p^s}, \ldots , \d_n^{p^s})$ is the $n$-tuple of
commuting inner derivations of $\CDPn$  such that $\d_i^{p^s}(\D_n)
\subseteq \D_n$ for all $i$. Let  $\ID (\d^{p^s}, \D_n)$ be the set
of all iterative $\d^{p^s}$-descents in the algebra $\D_n$ of
exponent $(\infty , \ldots , \infty )$. It follows from the defining
relations (\ref{DPndef}) for the $K$-algebra $\CDPn$  that the map
\begin{equation}\label{FIDs}
\FrobsCDPnDn \ra \ID (\d^{p^s}, \D_n), \;\; G\mapsto \{
G(\der^{[\alpha]})\}_{\alpha \in \Nn},
\end{equation}
is a {\em bijection} where the inverse map is given by the rule
$$\{ y^{[\alpha]}\}\mapsto G:\der^{[\alpha]}\mapsto y^{[\alpha]},
\;\; \alpha \in \Nn .$$ For each $i=1, \ldots , n$, let us introduce
the $K$-linear map 
\begin{equation}\label{sgid}
\int_i^*:\D_n\ra \D_n, \;\; \der^{[\alpha]}\mapsto \der^{[\alpha
+e_i]}, \;\; \alpha \in \Nn ,
\end{equation}
where $e_1, \ldots , e_n$ is the standard $\Z$-basis of
$\Z^n=\oplus_{i=1}^n\Z e_i$. The map    $\int_i^*$ is a right
inverse of the map  $\d_i$, $\d_i\int_i^*= \id_{\D_n}$. We call
the map $\int_i^*$ the $i$'th {\em   dual   integration} on $\D_n$
by analogy with the usual  $i$'th integration $\int_i$ on $P_n$
 when  ${\rm char}(K)=0$,  which is the $K$-{\em linear} map
$$\int_i:P_n\ra P_n, \;\; x^{[\alpha]}\mapsto  x^{[\alpha +e_i]},
\;\; x^{[\alpha]}:= \prod_{i=1}^n \frac{x^{\alpha_i}}{\alpha_i!},
\;\; \alpha \in \Nn .$$ The map $\int_i$ is a right inverse to
$\der_i:= \frac{\der}{\der x_i}$ in $P_n$ as $\der_i\int_i =
\id_{P_n}$.

For each $i\neq j$, the maps $\int_i^*$ and $\d_j$ commute. The
commutator $[\d_i, \int_i^*]=  1-\int_i^*\d_i$ is the projection
onto the subalgebra $\D_{n, \widehat{i}}:=  \t_{j\neq i}\D_1(j)$ in
the decomposition $\D_n = \oplus_{j\geq 0}\D_{n,
\widehat{i}}\derij$.

$\noindent $

{\it Definition}. For each natural number $s$ such that $s\geq 1$,
let $\CM_{n,s}$ be the set of all $n\times \N$ matrices
$u=(u_{ij})$, $i=1, \ldots , n$, $j\in \N$, such that
$$u_{ij}\in  \oplus_{0\neq \alpha \in
C_s^n}K\der^{[\alpha]}\subset \oplus_{\alpha \in
C_s^n}K\der^{[\alpha]}= \D_n\cap \cap_{i=1}^s \ker (\ad
(x_i^{p^s}))\;\;\; {\rm (Lemma \; \ref{a28Jun7}.(6))}$$
 where
$ C_s^n:= \{ \alpha \in \Nn \, | \, {\rm all}\;\; \alpha_i<p^s\}
$. The set $C_s^n$ is the $n$-{\em dimensional discrete cube} of
size $p^s$. It contains $p^{ns}$ elements. An element $a\in \D_n$
satisfies the condition that $a^p=0$ iff $a\in
\D_{n,+}:=\oplus_{0\neq \alpha \in \N^n}K\der^{[\alpha]}$. Since
$u_{ij}\in \D_{n,+}$, we have $u_{ij}^p=0$.

$\noindent $

The next theorem gives explicitly all the elements of the set $\ID
(\d^{p^s}, \D_n)$.

\begin{theorem}\label{M5Jul7}
Let $K$ be a field of characteristic $p>0$ and $s\geq 1$. Then   the
map
$$ \CM_{n,s}\ra \ID (\d^{p^s}, \D_n), \;\;  u=(u_{ij})\mapsto  y_u:= \{
y_u^{[\alpha]} \}_{\alpha \in \N^n},$$ is a    bijection where
$y_u^{[\alpha]}:=\prod_{i=1}^n  y_{u,i}^{[\alpha_i]}$,
$y_{u,i}^{[j]}:= \prod_{k\geq 0} \frac{y_{u,i}^{[p^k]j_k}}{j_k!}$
(for each $j\in \N$ written $p$-adically as $j=  \sum j_kp^k$),
$y_{u,i}^{[1]}:=\derips +u_{i0}$, and then  recursively, for each
$k\geq 1$,
$$ y_{u,i}^{[p^k]}:=
u_{ik}+\int_i^{*p^s}\prod_{l=0}^{k-1}\frac{y_{u,i}^{[p^l](p-1)}}{(p-1)!}.$$
The map $\{ z^{[\alpha]}\}_{\alpha \in \N^n}\mapsto u=(u_{ik})$,
given by the rule
$$u_{ik}:= z_i^{[p^k]}-\int_i^{*p^s}z_i^{[p^k-1]}, $$
is the inverse map of the  map $u\mapsto y_u$.
\end{theorem}

{\it Remark}. Note that $y_{u,i}^{[p^k]}:=
u_{ik}+\int_i^{*p^s}y_{u,i}^{[p^k-1]}$, $k\geq 1$.

$\noindent $

{\it Proof}.  First, we prove that the map $u\mapsto y_u$ is
well-defined, i.e. $y_u\in \ID (\d^{p^s}, \D_n)$. The ring $\D_n$ is
commutative. By Corollary \ref{f14Nov7}, $y_u\in \ID (\d^{p^s},
\D_n)$  iff, for each each $i,j=1, \ldots , n$ and $k\in \N$,
$$ y_{u,j}^{[p^k]p}=0\;\; {\rm and}\;\; \d_i^{p^s}(y_{u,j}^{[p^k]})= \d_{i,j} y_{u,j}^{[p^k-1]},$$
where $\d_{i,j}$ is the Kronecker delta. By the very definition,
the elements $y_{u,i}^{[p^k]}$  belong to the set $\D_{n,+}$, and
so the first type of  equalities hold, i.e.  $
y_{u,j}^{[p^k]p}=0$. To verify that  the second type of equalities
hold we use induction on $k$. The case $k=0$ is trivially true.
Suppose that $k\geq 1$ and that the equalities hold for  all $k'$
such that $k'<k$. Then
$$ \d_i^{p^s}(y_{u,i}^{[p^k]})=
\d_i^{p^s}(u_{ik}+\int_i^{*p^s}y_{u,i}^{[p^k-1]})=
0+\d_i^{p^s}\int_i^{*p^s}y_{u,i}^{[p^k-1]}=y_{u,i}^{[p^k-1]},$$ and,
for $i\neq j$,
$$ \d_j^{p^s}(y_{u,i}^{[p^k]})=\d_j^{p^s}(u_{ik}+\int_i^{*p^s}y_{u,i}^{[p^k-1]})=
0+\int_i^{*p^s}\d_j^{p^s}y_{u,i}^{[p^k-1]}=\int_i^{*p^s}0=0.$$ By
induction,  the equalities hold for all $k\in \N$. This shows that
the map $u\mapsto y_u$ is well-defined.

By the very definition, the map $u\mapsto y_u$ is injective. To
prove that it is surjective we have  to show that each element
$y=\{ y^{[\alpha]}\}\in \ID (\d^{p^s}, \D_n)$ is equal to $y_u$
for some $u=(u_{ik})\in \CM_{n,s}$. Let us set
$u_{ik}:=y_i^{[p^k]}-\int_i^{*p^s}y_i^{[p^k-1]}$. Clearly,
$u_{ik}^p=0$, $\d_i^{p^s}(u_{ik})= y_i^{[p^k-1]}-y_i^{[p^k-1]}=0$
and, for each $j\neq i$,
$$\d_j^{p^s}(u_{ik})=0-\int_i^{*p^s}\d_j^{p^s}y_i^{[p^k-1]}=0.$$
Therefore, $u_{ik}\in \sum_{0\neq \alpha \in
C_s^n}K\der^{[\alpha]}$. By the very definition, $y=y_u$. This
proves that the map $u\mapsto y_u$ is a bijection and $\{
z^{[\alpha]}\}\mapsto  u  = (u_{ik})$ is its inverse map. $\Box $

$\noindent $

As a consequence of Theorem \ref{M5Jul7}, we can find explicitly the
sets $\FrobsCDPnDn$, $s\geq 1$.

\begin{theorem}\label{cM5Jul7}
Let $K$ be a field of characteristic $p>0$ and $s\geq 1$. Then the
map
$$ \CM_{n,s}\ra \FrobsCDPnDn, \;\; u\mapsto
G_u:\der^{[\alpha]}\mapsto y_u^{[\alpha]}, \;\; \alpha \in \Nn, $$
is a bijection where the elements $y_u^{[\alpha]}$ are defined in
Theorem \ref{M5Jul7}, and the map
$$G\mapsto (u_{ik}(G)), \;\; u_{ik}(G):= G(\deripk ) -
\int_i^{*p^s}G(\der_i^{[p^k-1]}),$$ is its inverse  map.
\end{theorem}

{\it Proof}. The theorem    follows directly from Theorem
\ref{M5Jul7} and (\ref{FIDs}). $\Box $

$\noindent $

{\it Example}. Let $n=s=1$, $u=(\l \der , 0,0, \ldots )$, $\l \in
K$,  and $F':=G_u$. Then $F'(\der )= \der^{[p]}+\l \der $ and
\begin{eqnarray*}
  F'(\der^{[p]})&=& \int^{*p}\frac{(\der^{[p]}+\l\der)^{p-1}}{(p-1)!}=
 \frac{1}{(p-1)!}\int^{*p}\sum_{j=0}^{p-1}{p-1\choose j}\der^{[p]j}(\l \der )^{p-1-j}  \\
 &=&\int^{*p}\sum_{j=0}^{p-1}\l^{p-1-j}\der^{[jp]}\der^{[p-1-j]}=
 \int^{*p}\sum_{j=0}^{p-1}\l^{p-1-j}\der^{[p-1-j+jp]}\\
 &=& \sum_{j=0}^{p-1}\l^{p-1-j}\der^{[p-1-j+(j+1)p]}.
\end{eqnarray*}

The ring $\CDPn$ is simple, so each Frobenius $F'\in \FrobCDPn$ is
automatically a monomorphism.

\begin{corollary}\label{7Apr8}
Each Frobenius $F'\in \FrobCDPn$ is not an automorphism of the
ring $\CDPn$. Moreover, the ring $\CDPn$ is left and right free
finitely generated $KF'(\CDPn )$-module of rank $p^{2n}$.
\end{corollary}

{\it Proof}. By Theorem \ref{18Mar8}.(1), we may assume that
$F'\in \FrobCDPnPnDn$. Note that $\CDPn = P_n\t \D_n$,
$F'(P_n)\subseteq P_n$ and $F' (\D_n) \subseteq \D_n$. The
commutative algebra $P_n$ is a free $KF'(P_n)$-module of rank
$p^n$. By Theorem \ref{M5Jul7}, the commutative algebra $\D_n$ is
a free $KF'(\D_n)$-module of rank $p^n$. Then it is easy to deduce
 (using Theorem \ref{M5Jul7}) that
 $$ \CDPn = \oplus_{\alpha , \beta \in C^n_1}KF'(\CDPn )x^\alpha
 \der^{[\beta ] }=\oplus_{\alpha , \beta \in C^n_1}x^\alpha
 \der^{[\beta ] }KF'(\CDPn ),$$
 and the results follow. $\Box $


\section{Appendix: Some technical results}

In this section, some obvious technical results on the ring
$\CDPn$ are collected. Proofs are included for reader's
convenience. Let $\D_1(i):=\oplus_{j\geq 0} K\derij$.

\begin{lemma}\label{a28Jun7}
Let $K$ be a field of characteristic $p>0$; $k,k_i\in \N$, and
$\ker(\cdot ) := \ker_{\CD (P_n)}(\cdot )$. Then
\begin{enumerate}
\item  $\ker \, \ad (\deripk )= \CD (K[x_1, \ldots , \hx_i, \ldots
, x_n])\t K[x_i^{p^{k+1}}]\t \D_1(i)$.  \item $\cap_{i=1}^n \ker
\, \ad (\der_i^{[p^{k_i}]})= \D_n\t K[x_1^{p^{k_1+1}}, \ldots ,
x_n^{p^{k_n+1}}]$. \item $\cap_{k\geq 0} \ker \, \ad (\deripk ) =
\CD (K[x_1, \ldots , \hx_i, \ldots , x_n])\t \D_1(i)$. \item
$\cap_{i=1}^n\cap_{k\geq 0} \ker \, \ad (\deripk ) =\D_n$. \item
$\ker \, \ad (x_i^{p^k}) =\CD (K[x_1, \ldots , \hx_i, \ldots ,
x_n])\t K[x_i]\t K[ \der_i, \der_i^{[p]}, \ldots ,
\der_i^{[p^{k-1}]}]$. \item $\cap_{i=1}^n \ker \, \ad
(x_i^{p^{k_i}})= P_n \t \t_{i=1}^n K[ \der_i, \der_i^{[p]}, \ldots
, \der_i^{[p^{k_i-1}]}]$. \item $\cap_{i=1}^n \ker \, \ad (x_i) =
P_n$.\item For each natural number $k\geq 1$, the centralizer
$C(F^k_x(\CD (P_n)), \CD (P_n))$ of the ring $F^k_x(\CD (P_n))$ in
$\CDPn$ is equal to the direct sum $\oplus_{\alpha} K\der^{[\alpha
]} $ where $\alpha = (\alpha_1, \ldots , \alpha_n)\in \Nn$, $0\leq
\alpha_1<p^k, \ldots , 0\leq \alpha_n<p^k$.

\end{enumerate}
\end{lemma}

{\it Proof}. 1. The RHS of the equality in question is a subset of
the LHS. The opposite inclusion follows at once from the following
equality: for each natural number $j$ written $p$-adically as $j =
\sum_{s\geq l } j_sp^s$ with $j_l\neq 0$ and such that $l\leq k$,
there is the equality  
\begin{equation}\label{cpkdi}
[\deripk , x_i^j] = (j_l\der_i^{[p^k-p^l]}x_i^{(j_l-1)p^l} +
\cdots ) x_i^{\sum_{s\geq k+1} j_sp^s}
\end{equation}
where the three dots mean an element of the set $\sum_{0\leq t
<p^k-p^l} \der_i^{[t]} K[x_i]$. In more detail,
\begin{eqnarray*}
 [\deripk , x_i^j]& =& [\deripk , \prod_{s\geq l} x_i^{j_sp^s}]= [\deripk , \prod_{s=l}^k x_i^{j_sp^s}]\cdot \prod_{t\geq k+1} x_i^{j_tp^t}\\
 &=& ([\deripk , x_i^{j_lp^l}]\prod_{\nu = l+1}^k x_i^{j_\nu p^\nu
 } +x_i^{j_lp^l}[ \deripk ,
 x_i^{j_{l+1}p^{l+1}}] \prod_{\nu = l+2}^k x_i^{j_\nu p^\nu
 }+\cdots \\
 & & + x_i^{j_lp^l}\cdots x_i^{j_{\mu -1}p^{\mu -1}}[\deripk ,
 x_i^{j_\mu p^\mu}]\prod_{\nu = \mu +1}^k x_i^{j_\nu p^\nu }
 +\cdots )\cdot  \prod_{t\geq k+1} x_i^{j_tp^t}.
\end{eqnarray*}
Since
\begin{eqnarray*}
 [\deripk , x_i^{j_\mu p^\mu}]& =& F_x^\mu ([\der_i^{[p^{k-\mu}]},
 x_i^{j_\mu}])=F_x^\mu (j_\mu \der_i^{[p^{k-\mu }-1]}x_i^{j_\mu
 -1}+a_\mu)\\
 &=& j_\mu \der_i^{[p^k-p^\mu]}x_i^{(j_\mu -1) p^\mu }
 +F_x^\mu (a_\mu )
\end{eqnarray*}
for some element $a_\mu \in \sum_{r<p^{k-\mu}-1}
\der_i^{[r]}K[x_i]$, the equality (\ref{cpkdi}) is obvious.

2 and 3. Statements 2 and 3 follow from statement 1.

4. Statement 4 follows from statement 3.

5. The RHS of the equality in statement 5  is a subset of the LHS.
The opposite inclusion follows directly from the following
equality: for each natural number $j$ written $p$-adically as $j =
\sum_{0\leq s \leq l} j_sp^s$ with $j_l\neq 0$ and such that
$l\geq k$ (i.e. $j\geq p^k$), 
\begin{equation}\label{1cpkdi}
[\derij , x_i^{p^k}]= \der_i^{[j-p^k]}.
\end{equation}
In more detail, let $a:= \der_i^{[\sum_{s=0}^{k-1}j_sp^s]}$, $b:=
\der_i^{[\sum_{s=k}^lj_sp^s]}= F_x^k (c)$ and  $c:=
\der_i^{[\sum_{s=k}^l j_sp^{s-k}]}$. Then
\begin{eqnarray*}
 [\derij , x_i^{p^k}]& =& [ab, x_i^{p^k}]= a[F_x^k(c) , F_x^k( x_i)]=
 aF_x^k([c,x_i])= a F_x^k(\der_i^{[\sum_{s=k}^lj_sp^{s-k}-1]} )\\
 &=&
 a\der_i^{[\sum_{s=k}^lj_sp^s-p^k]}=\der_i^{[j-p^k]}.
\end{eqnarray*}

6. Statement 6 follows from statement 5.

7. Statement 7 is a particular case of statement 6.

8. Note that the centralizer $C:=C(F_x^k(\CD (P_n)), \CD (P_n))$
is the centralizer of the set that is the image under $F_x^k$ of
the set of canonical generators for $\CDPn$. Let $C'$ be the sum
$\sum_\alpha K\der^{[\alpha ]}$ in statement 8 (we have to show
that $C=C'$). By statements 1 and 5,
$$C_1:= \cap_{i=1}^n \ker (x_i^{p^k})\cap \cap_{i=1}^n \ker ( \ad
\, \deripk )= KF_x^{k+1} (P_n) \t C'.$$ Now, $C:= C_1\cap
\cap_{i=1}^n \cap_{l\geq k+1} \ker (\ad \, \der_i^{[p^l]}) = C'$.
 $\Box $

\begin{lemma}\label{b28Jun7}
Let $K$ be a field of characteristic $p>0$; $k,k_i\in \N$, and
$\ker(\cdot ) := \ker_{P_n}(\cdot )$. Then
\begin{enumerate}
\item  $\ker (\deripk )= K[x_1, \ldots , \hx_i, \ldots , x_n]\t
K[x_i^{p^{k+1}}]$.
 \item
$\cap_{i=1}^n \ker  (\der_i^{[p^{k_i}]} ) = K [ x_1^{p^{k_1+1}},
\ldots , x_n^{p^{k_n+1}}]$.
 \item
$\cap_{k\geq 0} \ker (\deripk ) =K[x_1, \ldots , \hx_i, \ldots ,
x_n]$. \item $\cap_{i=1}^n\cap_{k\geq 0} \ker (\deripk ) = K$.
\end{enumerate}
\end{lemma}

{\it Proof}. 1. The RHS of the equality in question is a subset of
the LHS. The opposite inclusion follows the  equality
$$ \deripk *x_i^j = {j\choose p^k} x_i^{j-p^k}.$$
Note that ${j\choose p^k}=0$ iff $j_k=0$   in the  $p$-adic   sum
$j = \sum   j_sp^s$.

2 and 3. Statements 2 and 3 follow  from  statement 1.

4. Statement 4  follows from statement  3. $\Box$

$${\bf Acknowledgements}$$

A part of the paper was written during my two months visit to the
IHES in 2007, the hospitality and support of the IHES is greatly
acknowledged.  The author would like to thank P. Cartier, O.
Gabber and M. Kontsevich for several discussions on various
aspects of rings of differential operators in prime characteristic
and their modules and for their comments.

Department of Pure Mathematics

University of Sheffield

Hicks Building

Sheffield S3 7RH

UK

email: v.bavula@sheffield.ac.uk







%


\begin{thebibliography}{99}

\bibitem{BCW}    H. Bass, E. H. Connell and D. Wright, The Jacobian Conjecture:
reduction of degree and formal expansion of the inverse, {\it
Bull. Amer. Math. Soc. (New Series)},
 {\bf 7} (1982), 287--330.



\bibitem{Bav-simderharp} V. V. Bavula, Simple derivations of differentiably simple Noetherian commutative rings in prime
characteristic,  {\em Trans. Amer. Math. Soc.},  {\bf 360} (2008),
no. 8, 4007-4027. (Arxiv:math.RA/0602632).

\bibitem{Bav-cpinv} V. V. Bavula,   The inversion formulae for automorphisms of polynomial
algebras and differential operators in prime characteristic, {\em
J. Pure Appl. Algebra}, {\bf 212} (2008), 2320--2337.
(Arxiv:math.RA/0604477).

\bibitem{Bav-ProcAMS2009} V. V. Bavula, The group of order preserving automorphisms of the ring of
differential operators on a Laurent polynomial algebra in prime
characteristic, {\em Proc. Amer. Math. Soc.}, {\bf 137} (2009),
1891-1898, (Arxiv:math.RA/0806.1038).

\bibitem{JC-DP} V. V.  Bavula, The ${\rm Jacobian \; Conjecture}_{2n}$ implies the ${\rm Dixmier
\;Problem}_n$, ArXiv:math.RA/0512250.


\bibitem{Bel-Kon05JCDP} A. Belov-Kanel and M. Kontsevich, The
Jacobian conjecture is stably equivalent to the Dixmier
Conjecture, {\em Mosc. Math. J.},  {\bf 7} (2007), no. 2, 209--218
 (arXiv:math. RA/0512171).


\bibitem{Dix-1968} J. Dixmier, Sur les algebres de Weyl. {\em   Bull. Soc. Math. France}, {\bf  96}  1968, 209--242.

\bibitem{Tsuchi05} Y. Tsuchimoto, Endomorphisms of Weyl algebra and $p$-curvatures. {\em Osaka J. Math.},  {\bf 42} (2005), no. 2, 435--452.


\end{thebibliography}
\end{document}